\documentclass[12pt]{amsart}
\usepackage{a4wide}
\usepackage{amssymb}
\usepackage{amsthm}
\usepackage{amsmath}
\usepackage{amscd}

\usepackage{verbatim}
\usepackage[all]{xy}

\numberwithin{equation}{section}

\theoremstyle{plain}
\newtheorem{theorem}{Theorem}[section]
\newtheorem{corollary}[theorem]{Corollary}
\newtheorem{lemma}[theorem]{Lemma}
\newtheorem{proposition}[theorem]{Proposition}

\newtheorem{conjecture}[theorem]{Conjecture}
\theoremstyle{definition}

\newtheorem{remark}[theorem]{Remark}

\theoremstyle{remark}

\newcommand{\OO}{\mathcal O}
\newcommand{\A}{\mathbb{A}}
\newcommand{\R}{\mathbb{R}}
\newcommand{\G}{\mathbb{G}}
\newcommand{\Q}{\mathbb{Q}}
\newcommand{\Z}{\mathbb{Z}}

\newcommand{\C}{\mathbb{C}}

\renewcommand{\H}{\mathbb{H}}

\newcommand{\zxz}[4]{\begin{pmatrix} #1 & #2 \\ #3 & #4 \end{pmatrix}}

\newcommand{\kzxz}[4]{\left(\begin{smallmatrix} #1 & #2 \\ #3 & #4\end{smallmatrix}\right) }

\font\cute=cmitt10 at 12pt 
\newcommand{\kay}{{\text{\cute k}}}

\newcommand{\frakp}{\mathfrak p}

\newcommand{\norm}{\operatorname{N}}

\newcommand{\vol}{\operatorname{vol}}
\newcommand{\tr}{\operatorname{tr}}

\newcommand{\Cl}{\operatorname{Cl}}

\newcommand{\GSpin}{\operatorname{GSpin}}

\newcommand{\Pet}{\text{\rm Pet}}
\newcommand{\Gr}{\operatorname{Gr}}
\newcommand{\GL}{\operatorname{GL}}
\newcommand{\SO}{\operatorname{SO}}

\newcommand{\Res}{\operatorname{Res}}

\newcommand{\z}{\operatorname{Z}}

\newcommand{\cha}{\operatorname{{Char}}}

\newcommand{\ord}{\operatorname{ord}}
\newcommand{\Gspin}{\operatorname{GSpin}}

\newcommand{\ff}{\hbox{if }}
\newcommand{\SL}{\operatorname{SL}}
\newcommand{\Diff}{\operatorname{Diff}}

\newcommand{\kk}{\kay}

\newcommand{\ph}{\phi}

\newcommand{\Ind}{\operatorname{Ind}}

\newcommand{\half}{\frac{1}{2}}

\newcommand{\CM}{\operatorname{CM}}

\begin{document}

\title[The lambda invariants at CM points]{The lambda invariants at CM points}

\author[ Tonghai Yang, Hongbo Yin, Peng Yu]{Tonghai Yang, Hongbo Yin, Peng Yu}

\address{Department of Mathematics, University of Wisconsin Madison, Van Vleck Hall, Madison, WI 53706, USA}
\email{thyang@math.wisc.edu}
\address{School of Mathematics, Shandong University, Jinan 250100,
P.R.China}
\email{yhb2004@mail.sdu.edu.cn}
\address{Department of Mathematics, University of Wisconsin Madison, Van Vleck Hall, Madison, WI 53706, USA}
\email{yu@math.wisc.edu}

\subjclass[2000]{14G35, 14G40, 11G18, 11F27}

\thanks{The first author is partially supported by a NSF grant DMS-1762289. The second author is partially supported by NSFC-11701548.}


\begin{abstract} In the paper, we show that $\lambda(z_1) -\lambda(z_2)$, $\lambda(z_1)$ and $1-\lambda(z_1)$ are all Borcherds products in $X(2) \times X(2)$. We then use    the big CM value formula of Bruinier, Kudla, and Yang to give  explicit factorization formulas for the norms of  $\lambda(\frac{d+\sqrt d}2)$,  $1-\lambda(\frac{d+\sqrt d}2)$, and $\lambda(\frac{d_1+\sqrt{d_1}}2) -\lambda(\frac{d_2+\sqrt{d_2}}2)$, with  the latter  under the condition $(d_1, d_2)=1$. Finally, we use these results to show that  $\lambda(\frac{d+\sqrt d}2)$ is always an algebraic integer and  can be easily used to construct units in the ray class field of $\Q(\sqrt{d})$ of modulus $2$.    In the process, we also give explicit formulas for a whole family of local Whittaker functions, which are of independent interest.
\end{abstract}

\maketitle

\setcounter{tocdepth}{1}
\tableofcontents

\section{Introduction}

It is well-known that the CM value $j(\tau)$ of the $j$-invariant at a CM point $\tau$ is always integral. What about the $\lambda$-invariants on  the  modular curve $X(2)$? Actually, there are six of them, satisfying the equation over $\Q(j)$
\begin{equation} \label{eq:generic}
f(\lambda, j) =  (1-\lambda + \lambda^2)^3 -\frac{ j}{256}  \lambda^2 (1-\lambda)^2 =0.
\end{equation}
The group $\SL_2(\Z)/\Gamma(2) =S_3$ acts on these roots simply transitively. The standard choice is
\begin{equation} \label{eq:lambda}
\lambda(\tau)=-\frac{1}{16}q^{-\frac{1}{2}}\prod_{n\geq1}(\frac{1-q^{n-\frac{1}{2}}}{1+q^n})^8, \quad q=e(\tau) =e^{2 \pi i \tau}.
\end{equation}
The other five roots of (\ref{eq:generic})  are $1-\lambda$, $\frac{1}{\lambda}$, $\frac{1}{1-\lambda}$, $\frac{\lambda-1}{\lambda}$ and $\frac{\lambda}{\lambda-1}$, all with   similar product formulas. For an imaginary quadratic field $\kay=\Q(\sqrt d)$ and a positive integer $N$, let $\CM(\kay, N)$ be the set of CM  points $\tau \in  Y(N)=\Gamma(N) \backslash \H$ such that the associated elliptic  curve $E_\tau=\C/(\Z +\Z\tau)$ has CM by $\mathcal O_\kay$. For every CM point $\tau \in \CM(\kay, 2)$, (\ref{eq:generic}) implies that $\lambda(\tau)$ is integral away from $2$. The following examples
\begin{align*}
\lambda(\frac{-7 +\sqrt{-7}}2) &=\frac{1+3\sqrt{-7}}{2} \hbox{ and } \,
 \lambda(\frac{-7 +\sqrt{-7}}4) =\frac{1+3\sqrt{-7}}{32}
\end{align*}
show that whether $\lambda(\tau)$ is integral  might not be a trivial question. Our first main result indicates that the standard choice of $\lambda(\tau)$ might be the best choice for the CM point $\frac{d+\sqrt d}2$.

\begin{theorem}  \label{theo:integral} Let $d<0$ be a fundamental discriminant of an  imaginary quadratic field $\kay=\Q(\sqrt d)$. Then
$\lambda_0 =\lambda(\frac{d+\sqrt d}2)$ is an algebraic integer.  Moreover,
\begin{enumerate}
\item  When $d \equiv  5 \pmod 8$, both $\lambda_0$  and $1-\lambda_0$ are units.

\item When $d \equiv 4 \pmod 8$, $\lambda_0$ is a unit.

\item  When  $d \equiv  0 \pmod 8$, $1-\lambda_0$ is a unit.

\item When  $d\equiv 1 \pmod 8$, $\frac{16}{\lambda_0(1-\lambda_0)}$ is a unit.
\end{enumerate}
\end{theorem}

The case $d \equiv 5 \pmod 8$ is relatively straightforward as one can prove that $[\Q(\lambda_0) :\Q(j_0)]=6$ with $j_0 =j(\frac{d+\sqrt d}2)$. The other cases are quite subtle. When $d \not\equiv 5 \pmod 8$, $\lambda_0$ is of degree $2$ over $\Q(j_0)$ and its Galois conjugate may be $1-\lambda_0$, $\frac{1}{\lambda_0}$, or $\frac{\lambda_0}{\lambda_0-1}$, but it is highly non-trivial to determine which one is its Galois conjugate in each individual case. The key technical result to prove the above theorem is

\begin{proposition}  \label{prop:GaloisConjugate} The Galois conjugate of $\lambda_0$ over $\Q(j_0)$ is  $1-\lambda_0$, $\frac{1}{\lambda_0}$ or $\frac{\lambda_0}{\lambda_0-1}$ according to $d\equiv 1, 4$ or $0 \pmod 8$.
\end{proposition}
  The way we prove this proposition is to use  the following norm formulas. Notice that Galois conjugates have to have the same norm.

\begin{proposition}  \label{prop:norm} We have the following norm formulas  for $\lambda_0 =\lambda(\frac{d+\sqrt{d}}2)$ with  $ d <-4$. Here $h$ is the class number of $\kay=\Q(\sqrt d)$.
\begin{table}[h]
\caption{\label{table1}$|\mathrm{N}(\lambda)|$}
\begin{tabular}{|c|c|c|c|c|}
\hline \rule[-3mm]{0mm}{8mm}
$d\pmod 8$ & $0$ & $ 1$ &$4$ &$5$\\
\hline $|\mathrm{N}(\lambda_0)|$    &$2^{2h}$    &$2^{4h}$   &$1$    &$1$\\
\hline $|\mathrm{N}(1-\lambda_0)|$  &$1$      &$2^{4h}$   &$2^{2h}$  &$1$\\
\hline $|\mathrm{N}(1/\lambda_0)|$  &$2^{-2h}$ &$2^{-4h}$  &$1$     &$1$\\
\hline $|\mathrm{N}(\lambda_0/(\lambda_0-1))|$  &$2^{2h}$    &$1$        &$2^{-2h}$  &$1$\\
\hline
\end{tabular}
\end{table}
\end{proposition}

This proposition is a bi-product of the original purpose of this paper: a factorization formula for $|\norm( \lambda(\frac{d_1+\sqrt{d_1}}2)-\lambda(\frac{d_2+\sqrt{d_2}}2))|$ similar to the beautiful factorization formula of singular moduli of Gross and Zagier (\cite{GZSingular}).

\begin{theorem} \label{maintheorem} Let $\kay_j=\Q(\sqrt{d_j})$ with ring of integers $\OO_{d_j} = \Z[\tau_j]$ with $ \tau_j =\frac{d_j +\sqrt{d_j}}2$, $j=1, 2$. Let $w_j$ be the number of roots of unity in $\kay_j$, and let $h_j=h(\kay_j)$ be the class number of $\kay_j$.  Let $E=\Q(\sqrt{d_1}, \sqrt{d_2})$ and $F=\Q(\sqrt D)$ with $D=d_1d_2$. For an ideal $\mathfrak a$ of $F$, let $\rho_{E/F}(\mathfrak a)$ be the number of integral ideals of $E$ with relative norm $\mathfrak a$, and let $\chi=\chi_{E/F}$ be the quadratic Hecke character of $F$ associated to $E/F$. Let $c(t) = \frac{d_1^2+d_2^2 -(d_1+d_2)}4 -\frac{m+D}2 \in \Z$. Assume that $(d_1 , d_2) =1$.
Then

\begin{enumerate}
\item When $d_1 \equiv d_2 \equiv 5 \pmod 8$,

\begin{align*}
&\log|\norm(\lambda(\tau_1) -\lambda(\tau_2))|
\\
&= 6
 \sum_{ \substack{ t =\frac{m+\sqrt D}2  \in  \OO_F \\ |m| < \sqrt D\\ c(t)\in 2\Z }}
 \sum_{\chi_{E/F}(\mathfrak p)=-1}\left(\frac{1+ \ord_{\mathfrak p_t}( t\mathfrak p_t^{-2})}2+\frac{\ord_{\frakp_t}(\frakp)}{3}\right)
  \rho_{E/F}(t\mathfrak p \mathfrak p_t^{-2})\log\norm(\mathfrak p).
\end{align*}
Here $\mathfrak p_t$ is the prime of $F$ above $2$ with $t \in \mathfrak p_t$.

\item  When  $d_1 \equiv 5 \pmod 8$, and $d_2 \equiv0, 4 \pmod 8$,

\begin{align*}
&\log|\norm(\lambda(\tau_1) -\lambda(\tau_2))|=
\sum_{ \substack{ t =\frac{m+\sqrt D}2  \in  \OO_F \\ |m| < \sqrt D\\ c(t)\in 2\Z }}
\sum_{\chi_{E/F}(\mathfrak p)=-1} \frac{1+ \ord_{\mathfrak p}(t  \mathfrak)}2
  \rho_{E/F}(t\mathfrak p \mathfrak)\log\norm(\mathfrak p).
\end{align*}

\item When $d_1 \equiv 1 \pmod 8$ and $d_2 \equiv 5 \pmod 8$,

\begin{align*}
&\log|\norm(\lambda(\tau_1) -\lambda(\tau_2))|= 2
\sum_{ \substack{ t =\frac{m+\sqrt D}2  \in  \OO_F \\ |m| < \sqrt D\\ c(t)\in 2\Z }}
\sum_{\chi_{E/F}(\mathfrak p)=-1} \frac{1+ \ord_{\mathfrak p}(t  \mathfrak)}2
  \rho_{E/F}(t\mathfrak p \mathfrak)\log\norm(\mathfrak p).
\end{align*}

\item When $d_1 \equiv 1 \pmod 8$ and $d_2 \equiv 0,4 \pmod 8$,

\begin{align*}
&\log|\norm(\lambda(\tau_1) -\lambda(\tau_2))|
\\
&=\frac{1}2
\sum_{ \substack{ t =\frac{m+\sqrt D}2  \in  \OO_F \\ |m| < \sqrt D\\ c(t)\in 2\Z }}
\sum_{\chi_{E/F}(\mathfrak p)=-1} \frac{1+ \ord_{\mathfrak p}(t  \mathfrak)}2
  \rho_{E/F}(t\mathfrak p \mathfrak)\log\norm(\mathfrak p)+2^{3}\frac{h_1 h_2}{w_1  w_2}\log 2.
\end{align*}

\item When $d_1 \equiv 1 \pmod 8$ and $d_2 \equiv 1 \pmod 8$,

\begin{align*}
&\log|\norm(\lambda(\tau_1) -\lambda(\tau_2))|
\\
&=2
\sum_{ \substack{ t =\frac{m+\sqrt D}2  \in  \OO_F \\ |m| < \sqrt D\\ c(t)\in 2\Z }}
\sum_{\chi_{E/F}(\mathfrak p)=-1} \frac{1+ \ord_{\mathfrak p}(t \mathfrak p_t^{-2})}2
  \rho_{E/F}(t\mathfrak p  \mathfrak p_t^{-2})\log\norm(\mathfrak p)+2^{5}\frac{h_1 h_2}{w_1 w_2}\log 2.
\end{align*}
Here $\mathfrak p_t$ is the prime of $F$ above $2$ with $t \in \mathfrak p_t$.

\end{enumerate}
\end{theorem}

 We remark  that for each $t =\frac{m+\sqrt D}2 \in \OO_F$ in the above formula, the second sum $\sum_{\chi_{E/F}(\mathfrak p) =-1}$ has at most one non-zero   term because it is impossible to have two primes $\mathfrak p_1$ and $\mathfrak p_2$ of $F$ inert in $E/F$ with  both $\rho_{E/F}(\mathfrak a  \mathfrak p_1)$ and $\rho_{E/F}(\mathfrak a \mathfrak p_2) $ non-zero for any fractional ideal $\mathfrak a$. In \cite{YuPeng}, one of the authors (Yu) extends the factorization formula to Rosenhain invariants on Siegel threefold.

  \begin{corollary} \label{cor1.2} Let the notation be as above.
  \begin{enumerate}

  \item $\norm(\lambda(\frac{d_1+\sqrt{d_1}}2) -\lambda(\frac{d_2+\sqrt{d_2}}2))$ is an integer.

  \item If $p|\norm(\lambda(\frac{d_1+\sqrt{d_1}}2) -\lambda(\frac{d_2+\sqrt{d_2}}2))$, then $p \le  \frac{D}4$ and $p | \frac{D -m^2}4$ for some integer $0 \le m < \sqrt D$.


  \item When $E/\Q$ is unramified at $2$,  $ \norm(\lambda(\frac{d_1+\sqrt{d_1}}2) -\lambda(\frac{d_2+\sqrt{d_2}}2))$ is a square up to sign.

  \end{enumerate}

  \end{corollary}

Although  $\norm(j(\tau_1) -j(\tau_2))$ does not depend on the choice of $\tau_j \in \CM(\kay_j, 1)$, the same is NOT true for $\norm(\lambda(\tau_1) -\lambda(\tau_2))$. Let
$$
\mathcal N(d_1, d_2) =\{ |\norm(\lambda(\tau_1) -\lambda(\tau_2))|:\, \tau_j \in  \CM(\kay_j, 2)\}.
$$
Then a simple group theoretic argument using Lemma \ref{lem:ClassGroup} shows that $\mathcal N(d_1, d_2)$  has order at most  $9$. We find the following surprising result.

\begin{theorem}  \label{theo:Difference} Assume $d_1, d_2 <-4$.  Then  the order of $\mathcal N(d_1, d_2)$ is given by
$$
|\mathcal N(d_1,d_2)|=
\begin{cases}
1&\text{if }d_1\equiv d_2 \equiv 5 \pmod 8,\\
2&\text{if exactly  one of }d_1,d_2\equiv 5\pmod 8,\\
4&\text{if }d_1\equiv d_2  \equiv 1\pmod 8,\\
5&\text{otherwise}.
\end{cases}
$$
\end{theorem}

  The proofs of Theorem  \ref{maintheorem} and Proposition \ref{prop:norm} are similar to that in \cite{YYCM} although they are technically more subtle and raise some interesting technical questions. In Section \ref{sect:Lambda1}, we review and prove some basic facts about CM points in $X(2)$, $\lambda$-invariants, and ray class group actions on the CM cycles. After that we prove Theorem  \ref{theo:integral} and Propositions \ref{prop:GaloisConjugate} assuming Proposition \ref{prop:norm}. In Section \ref{sect:Lambda}, we treat
  the product $X(2) \times X(2)$ of modular curves as the Shimura variety associated to $L=M_2(\Z)$ with $Q(x) =2 \det x$, and find a weakly holomorphic modular form $f$ whose Borcherds product lifting is $\lambda(z_1) -\lambda(z_2)$ (Proposition \ref{prop6.2}). In the process, we also prove that $\lambda(z_1)$ and $1-\lambda(z_1)$ are  Borcherds products for some explicit weakly homomorphic forms (Corollary \ref{lambda}). Notice that the existence of $f$ is not  guaranteed by Bruinier's converse theorem (\cite{Bruinier14}).
In Section \ref{sect:BigCM}, we use the big CM value formula in \cite{BKY12} to compute the absolute value of the norm of $\lambda(\tau_1) -\lambda(\tau_2)$ and to prove Theorem  \ref{maintheorem} and Proposition \ref{prop:norm}, although the big CM cycle and Galois orbit are not always the same in this case. We then  use these results to prove Theorem \ref{theo:Difference}.

  The main technical part to prove Theorem  \ref{maintheorem} and Proposition \ref{prop:norm} is to compute the local Whittaker functions at primes dividing $2$, which needs more general explicit formulas for local Whittaker functions than exists in literature. In Section \ref{sect5}, we give an explicit formula for a whole family of local Whittaker functions which are basic and of independent interest. In the last section, we actually compute some examples of the factorization of  the norm $|\norm(\lambda(\tau_1) -\lambda(\tau_2))|$ using Theorem
  \ref{maintheorem}.  We check the formulas with computer computation results, and they match perfectly.

\section{The integrality of the CM values of the  Lambda function $\lambda$ } \label{sect:Lambda1}

Let $Y(2) =\Gamma(2) \backslash \H$ and $X(2) = Y(2) \cup \{0, 1, \infty\}$. Then $X(2)$ is the compactification of $Y(2)$ and is known to have a canonical model over $\Q$. Then the $\lambda$-invariant defined by (\ref{eq:lambda})
gives an isomorphism between $X(2)$ and $\mathbb P^1$. Moreover, we have the following commutative diagram

\begin{equation}
\xymatrix{
  X(2)\ar[r]^-{\lambda} \ar[d]^-{\pi} &\mathbb P^1
  \ar[d]^-{\pi'} \cr
  X(1)\ar[r]^-{j} &\mathbb P^1 \cr }
   \end{equation}
where $\pi$ is the natural projection and   $\pi'$ is given by
$$
\pi'(x) =\frac{256 (1-x + x^2)^3}{x^2 (1-x)^2}.
$$
So the minimal polynomial of $\lambda=\lambda(\tau)$  over $\Q(j(\tau))$ is generically given by (\ref{eq:generic}).
For any $\tau \in \H$, the preimage of $j(\tau)$ under $\pi$ is
\begin{align*}
&\left\{ \lambda(\tau),  1-\lambda(\tau), \frac{1}{\lambda(\tau)}, \frac{1}{1-\lambda(\tau)},  \frac{\lambda(\tau)-1}{\lambda(\tau)},  \frac{\lambda(\tau)}{\lambda(\tau)-1}\right\}
\\
 &=\{ \lambda\circ \gamma: \,  \gamma \in \SL_2(\Z)/\Gamma(2)\}
 \end{align*}
Let $\kay=\Q(\sqrt d)$ be an imaginary quadratic field with fundamental discriminant $d$ viewed as a subfield of $\C$ in the usual way. For $N=1, 2$, let $\CM(\kay, N)$ be the set of CM points in
$X(N)$ defined as in the introduction. They are  parameterized by  triples $(\mathfrak a, \alpha, \beta)$, where $\mathfrak a$ is a fractional ideal of $\kay$, and $(\alpha, \beta)$ is an ordered $\Z$-basis of $\mathfrak a$ such that
$\frac{\alpha}{\beta}\in \H$. Two such triples $(\mathfrak a_i, \alpha_i, \beta_i)$ are equivalent if there is $z\in \kay$ and $\gamma \in \Gamma(N)$ such that
$$
\mathfrak a_1 =z \mathfrak a_2, \quad \gamma (\alpha_1, \beta_1)^t = (z \alpha_2, z\beta_2)^t.
$$
The associated CM point in $X(N)$ is given by $\tau_{(\mathfrak a, \alpha, \beta)}=\frac{\alpha}{\beta}$. The ray class group $\Cl(\kay, N)$ acts on $\CM(\kay, N)$ as follows: For a fractional ideal $\mathfrak b$ prime to $N$,
$$
[\mathfrak b]*(\mathfrak a, \alpha, \beta) = (\mathfrak a \mathfrak b, \alpha_1, \beta_1)
$$
with
$$
\alpha_1 \equiv \alpha \pmod N,  \hbox{ and } \beta_1 \equiv \beta \pmod N.
$$
We refer to \cite{Yang16} for more detail. Our special case $N=1, 2$ simplifies the description a little.  Although it is well-known that $\Cl(\kay)=\Cl(\kay, 1)$ acts on $\CM(\kay)=\CM(\kay, 1)$ simply transitively, it is not true even for $\Cl(\kay, 2)$.   Let $h(\kay, N)$ be the ray class number of $\kay$ of modulus $N$, and $h(\kay)=h(\kay,  1)$.

\begin{lemma} \label{lem:ClassGroup} Let the notation be as above, and let $w=|\OO_\kay^\times|$ be the number of roots of unity in $\kay$. Then
\begin{enumerate}
\item  One has
$$
\frac{h(\kay, 2)}{h(\kay)} =\frac{2}{w} r(d) =\begin{cases}
 \frac{2}w &\ff  d \equiv 1\pmod 8,
 \\
 \frac{4}w &\ff d \equiv 0 \pmod 4,
 \\
 \frac{6}w  &\ff d \equiv 5 \pmod 8.
  \end{cases}
$$
We denote $r(d)= 1,2, 3$ according to $d \equiv 1 \pmod 8 , 0 \pmod 4$, or $ 5 \pmod 8  $.

\item The CM set $\CM(\kay, 2)$ has $6/r(d)$ orbits under the action of $\Cl(\kay, 2)$ when $d <-4$ and has 3 orbits when $d=-3, -4$.

\item  When $4|d$, $\lambda(\tau)$ is real when $\tau \in \CM(\kay, 2)$ and the minimal polynomial of $\lambda(\tau)$ over $\Q$ is of degree $h(\kay, 2)$.

\item When $ 4\nmid d$, $\lambda(\tau)$  is not real for $\tau\in \CM(\kay, 2)$, and the minimal polynomials of $\lambda(\tau)$ and $\lambda(-\bar\tau)=\overline{\lambda(\tau)}$ over $\Q$ are the same and of degree $2h(\kay, 2)$. In particular,  when $d \equiv 5 \pmod 8$ and $d \ne -3$, the minimal polynomial of $\lambda(\tau)$ is given by
    $$
    f_\lambda(x) =\prod_{z \in \CM(\kay, 2)} (x-\lambda(z)).
    $$
\end{enumerate}

\end{lemma}

\begin{proof}(sketch)  Claim (1) follows from the general relation between  ideal class groups and ray class groups, see for example \cite[Theorem 1.7, p. 146]{Milne}:
$$
\frac{h(\kay, \mathfrak m)}{h(\kay)} =\frac{\norm(\mathfrak m)}{[\OO_\kay^\times: U_{\mathfrak m}]}\prod_{\mathfrak p |\mathfrak m} (1-\norm(\mathfrak p)^{-1}),
$$
where $\mathfrak m$ is an integral ideal of $\kay$, $h(\kay, \mathfrak m)$ is the ray class number of $\kay$ of modulus $\mathfrak m$, and  $U_{\mathfrak m}$ is the subgroup of $\OO_\kay^\times$ consisting of units which are $1\pmod {\mathfrak m}$.

Since the projection map from $X(2)$ to $X(1)$ has degree $6$ and it does not ramify at the CM points unless $d=-3, -4$,  claim (2)  follows from (1) immediately.

For (3), we may assume $\tau=\frac{-d+\sqrt d}2=[\OO_\kay, \alpha, 1]$ with $\alpha= \frac{-d+\sqrt d}2$.  Then $-\bar\tau=\frac{d+\sqrt d}2 = [\OO_\kay, -\bar\alpha, 1]$. When  $4|d$, it is clear $-\bar\tau= d+\tau$ which is  equal to $\tau$ in $X(2)$. Now assume that $d$ is odd.  If $\tau$  and $-\bar\tau$ are in the same orbit, then there is an ideal $\mathfrak b$ such that
$$
[\OO_\kay, -\bar\alpha, 1]=[\mathfrak b]*[ \OO_\kay, \alpha, 1]=[\mathfrak b, \alpha_1, \beta_1]
$$
So $\mathfrak b=z\OO_\kay$  and there is $\gamma \in \Gamma(2)$ such that
$$
\gamma (\alpha_1, \beta_1)^t = (-z \bar\alpha, z)^t.
$$
Notice that $\alpha_1 \equiv \alpha \pmod 2$ and $\beta_1 \equiv 1 \pmod 2$ by definition. A simple calculation shows that $\alpha +\bar \alpha=-d \equiv 0 \pmod 2$, a contradiction.

(4)  follows from (3).
\end{proof}

The cases $d=-3, -4$  are special, so we record them as a separate lemma.
\begin{lemma} One has for $d=-3, -4$
\begin{enumerate}
\item \begin{align*}
h(\Q(\sqrt{-4}), 2)&=h(\Q(\sqrt{-4})=1,  \quad  \lambda(\sqrt{-1})=-1,
\\
h(\Q(\sqrt{-3}), 2)&=h(\Q(\sqrt{-3})=1,  \quad \lambda(\frac{-3+\sqrt{-3}}2)=\frac{1+\sqrt{-3}}2.
\end{align*}

\item $\CM(\kay, 2)$ has three $\Cl(\kay, 2)$-orbits, and $X(2)\rightarrow X(1)$ is ramified at the CM points $\sqrt{-1}$ and $\frac{-3+\sqrt{-3}}2$ with ramification index $2$.

\end{enumerate}

\end{lemma}

\begin{proposition}
\label{prop:MininalPoly}
Let $d < -4$, and let $\tau =\tau_{[\mathfrak a, \alpha, \beta]}$ be a CM point in $X(2)$. Then
the polynomial in $\lambda$ $f(\lambda, j(\tau))$ in  (\ref{eq:generic}) is irreducible over $\Q(j(\tau))$ when  $d \equiv 5 \pmod 8$ and
factorizes into product of three quadratic polynomials over $\Q(j(\tau))$ when $d\not\equiv 5\pmod 8$. Equivalently,
$$
[\Q (\lambda(\tau)) : \Q(j(\tau))]=\begin{cases}
   6 &\ff d \equiv 5\pmod 8,
   \\
    2 &\hbox{ otherwise}.
    \end{cases}
$$
Moreover, when $d \not\equiv 5 \pmod 8$, we have
$$
f(\lambda, j(\tau)) =f_1(\lambda) f_2(\lambda) f_2(\lambda)
$$
with
\begin{align*}
f_1(\lambda) &=\lambda^2 -\lambda +a,
\\
f_2(\lambda)&=\lambda^2-\frac{1}a\lambda+\frac{1}a,
\\
f_3(\lambda) &=\lambda^2-(2-\frac{1}a)\lambda+1,
\end{align*}
for some $a \in \Q(j(\tau))$.
\end{proposition}
\begin{proof}(sketch) When  $4\nmid d$,  $[\Q(\lambda(\tau)):\Q] =2 h(\kay, 2)$ while $[\Q(j(\tau)):\Q]=h$, and so
$$
[\Q(\lambda(\tau)): \Q(j(\tau))]= 2 h(\kay, 2)/h,
$$
which is $6$ or $2$ depending when  $d \equiv 5$  or  $1 \pmod 8$  by Lemma \ref{lem:ClassGroup}. Similarly when  $4|d$, $[\Q(\lambda(\tau)): \Q(j(\tau))]= h(\kay, 2)/h=2$.

Now assume $d \not\equiv 5 \pmod 8$. Let $F=\Q(j(\tau))$. Since $F(\lambda(\tau))/F$ is of degree $2$,  the Galois conjugate of $\lambda(\tau)$ is either $1 -\lambda(\tau)$, $\frac{1}{\lambda(\tau)}$ or $\frac{\lambda(\tau)}{\lambda(\tau)-1}$.

First assume that the Galois conjugate of $\lambda(\tau)$ is  $1 -\lambda(\tau)$, then $a =\lambda(\tau) (1-\lambda(\tau)) \in F$, and the minimal polynomial of $\lambda(\tau)$ is $f_1(\lambda)$. At the same time, the Galois conjugate of $\frac{1}{\lambda}$ becomes $\frac{1}{1-\lambda}$, and so the minimal polynomial of $\frac{1}{\lambda(\tau)}$ is
$f_2(\lambda)$. Finally, the Galois conjugate of $\frac{\lambda(\tau)}{\lambda(\tau)-1} $ becomes $\frac{\lambda(\tau) -1}{\lambda(\tau)}$, and their  minimal polynomial is  $f_3(\lambda)$. As the roots of $f(\lambda, j(\tau)$ are these six numbers, we see $f(\lambda))=f_1 f_2 f_3$ as claimed. The other cases are similar.

\end{proof}

For the rest of this section, we assume Proposition  \ref{prop:norm} and use it  to prove Proposition \ref{prop:GaloisConjugate} and Theorem \ref{theo:integral}.  Proposition \ref{prop:GaloisConjugate} follows from
Propositions \ref{prop:norm} and \ref{prop:MininalPoly} immediately as Galois conjugates have to have the same norm.
Let $H(\kay, N)$ be the ray class field of $\kay$ of modulus $N$.

{\bf Proof of Theorem \ref{theo:integral}}:  We first prove that $\lambda_0=\lambda(\frac{d+\sqrt d}2)$ is integral.
  As $j(\tau)$ is an algebraic integer, the equation
$ (\ref{eq:generic})$ shows that $\lambda(\tau) $ is integral at $\mathfrak P$ for any CM point $\tau \in  \CM(\kay, 2)$ and every prime ideal $\mathfrak  P$ of $\kay(\lambda(\tau))=H(\kay, 2)$ not above $2$.
So we just need to prove that $\lambda_0$ is integral at $\mathfrak p$ for every prime $\mathfrak p$ of the Hilbert class field $H=\kay(j_0)$  of $\kay$ above $2$ with $j_0=j(\frac{d+\sqrt d}2)$.

When  $d \equiv 5 \pmod 8$, $[\kay(\lambda_0): \kay(j_0)] =3$ by  Lemma \ref{lem:ClassGroup}. The Galois conjugates of $\lambda_0$ over $\kay(j_0)=H$ are $\lambda_0$, $\frac{\lambda_0-1}{\lambda_0}$, and $\frac{1}{1-\lambda_0}$. So $\norm_{\kay(\lambda_0)/\kay(j_0)}(\lambda_0) =-1$. On the other hand, $2\mathcal O_\kay$ splits completely in $\kay(j_0)$ and has then to be ramified in $\kay(\lambda_0)=H(\kay, 2)$. So for every prime ideal $\mathfrak p$ of $\kay(j_0)$ above $2$, $\kay(\lambda_0)_{\mathfrak p}/\kay(j_0)_{\mathfrak p}$ is a ramified local field extension of degree $3$. So $\norm_{\kay(\lambda_0)/\kay(j_0)}(\lambda_0) =-1$ implies that $\lambda_0$ is a unit in $\kay(\lambda_0)_{\mathfrak p}$. Therefore, $\lambda_0$ is integral and actually a unit in $\kay(\lambda_0)$. So all six $\lambda$-values are units as they are Galois conjugates over $\Q(j_0)$.

When  $d \equiv 4 \pmod 8$, let $\mathfrak p_2 $ be the prime of $\kay$ above $2$. Then $\mathfrak p_2$ is unramified in $H$ and has to be ramified in $\kay(\lambda_0)=H(\kay, 2)$ (as $H \ne H(\kay, 2)$). In particular, for every prime $\mathfrak p$ of $H$ above $2$,  $\mathfrak p = \mathfrak P^2$ is ramified in  $\kay(\lambda_0)$.
By Proposition \ref{prop:GaloisConjugate}, $\lambda_0$ and $1/{\lambda_0}$ are Galois conjugates over $H$,  and so  we have
$$
\norm_{H(\kay, 2)_{\mathfrak P}/H_{\mathfrak p}}\lambda_0 =1,
$$
which  implies that $\lambda_0$ is a unit in $H(\kay, 2)_{\mathfrak P}$ for every prime ideal $\mathfrak P$ of $H(\kay, 2)$. So $\lambda_0 $ is integral and actually a unit in $H(\kay, 2)$.

Next, we assume $d \equiv 0 \pmod 8$.  As  the Galois conjugate of $\lambda_0$ is $\frac{\lambda_0}{\lambda_0 -1}$ by Proposition \ref{prop:GaloisConjugate},  the Galois conjugate of $1-\lambda_0$ is  $\frac{1}{1-\lambda_0}$. So the same argument as  in the case $d \equiv 4 \pmod 8$ shows that $1-\lambda_0$ is integral and a unit in $H(\kay, 2)$.

Finally, we assume  $d \equiv 1 \pmod 8$. Let
$$a=\lambda_0(1-\lambda_0),~\omega_0=-\frac{16}{a}.$$
It is not hard to find that
$$j_0=\frac{(\omega_0+16)^3}{\omega_0},\text{ or }(\omega_0+16)^3-j_0\omega_0=0.$$
Hence, $\omega_0$ is integral by integrality of $j_0$.  By Proposition \ref{prop:GaloisConjugate}, $\lambda_0$ and $1-\lambda_0$ are in the same Galois orbit. By Proposition \ref{prop:norm},
$$|\mathrm N_{\kay(\lambda_0)/\Q}a|=|\mathrm N_{\kay(j_0)/\Q}\lambda_0|^2=|\mathrm N(\lambda_0)|^2=2^{8h}.$$
Here $h$ is the ideal class number of $\kay$.
Therefore,
$$|\mathrm N\omega_0|=\frac{16^{2h}}{|\mathrm N_{\kay(\lambda_0)/\Q}a|}=\frac{16^{2h}}{2^{8h}}=1.$$ So $\omega_0$ is a unit as claimed.  In  particular, $a =-\frac{16}{\omega_0}$ is integral, and thus $\lambda_0$ is integral with
$$
\lambda_0^2 - \lambda_0 -a =0.
$$
This proves the theorem.

We remark that $\omega_2(\tau) =-\frac{16}{\lambda(\tau) (1-\lambda(\tau))}$ is actually  a hauptmodul for $X_0(2)$. Its CM values are studied in \cite{YYCM}.

\begin{corollary} The following are true.

\begin{enumerate}

\item When  $d \equiv 5 \pmod 8$, one has $256|j(\frac{d+\sqrt d}2)$.

\item When  $d \equiv 1 \pmod 8$, $\hbox{gcd}(j(\frac{d+\sqrt d}2), 2) =1$.

\item When  $d \equiv 0 \pmod 4$, $2^{6h} \| \norm_{\Q(j_0)/\Q}(j_0) $ for $j_0 =j(\frac{d+\sqrt d}2)$ where $h$ is the class number of $\kay=\Q(\sqrt d)$. Moreover, $\frac{256}{j_0}$ is integral at $2$.

\end{enumerate}

\end{corollary}
\begin{proof} (1) follows from the minimal polynomial equation (\ref{eq:generic}) of $\lambda_0$ over $\Z(j_0)$ and the integrality of $\lambda_0$. Here we write again $\lambda_0=\lambda(\frac{d+\sqrt d}2)$.

For (2),  recall that  $\omega_0= -\frac{16}{a}$ is a unit with $a =\lambda_0 (1-\lambda_0) \in \Q(j_0)$. So
$\ord_{\mathfrak  p}(a) =4 \ord_{\mathfrak p} 2$ for every prime ideal $\mathfrak p$ of $\Q(j_0)$ above $2$. This implies
$$
\ord_{\mathfrak p}j_0 = \ord_{\mathfrak p} \frac{256 (1 - a)^3}{a^2} =0.
$$
So $j_0$ is a unit at $2$. Claim (3) can be proved similarly and we leave it to the reader.
\end{proof}

The explicit calculation suggests  that $\frac{j_0}{2^6}$ is  integral and is a unit at $2$.

\section{$\lambda(z_1)$, $1-\lambda(z_1)$ and  $\lambda(z_1) -\lambda(z_2)$ as   Borcherds liftings} \label{sect:Lambda}

\subsection{Brief review of Borcherds liftings} \label{sect:Borcherds}
Let $L=M_2(\Z)$ with the quadratic form $Q(x) =2 \det x$ and the induced pairing
\[(x,y)=Q(x+y)-Q(x)-Q(y)=2(y_1x_4+x_1y_4-y_2x_3-x_2y_3).\]
Let $V=M_2(\Q)$.Then
$$
H=\GSpin(V) =\{ g=(g_1, g_2)\in \GL_2 \times \GL_2: \,  \det g_1=\det g_2\}
$$
acting on $V$ via
$$
(g_1, g_2) x =g_1 x g_2^{-1}.
$$
One has the exact sequence
$$
1 \rightarrow  \mathbb G_m \rightarrow H \rightarrow \SO(V) \rightarrow 1.
$$

 Let
 \begin{equation}
 \mathcal L = \{ w \in V_\C:\,  (w, w) =0,  \quad (w, \bar w)<0\}.
 \end{equation}
 and let $\mathbb D$ be the Hermitian symmetric domain of oriented negative $2$-planes in $V_\R=V\otimes_\Q \R$. Then one has an isomorphism
 $$
 pr: \mathcal L/\C^\times \cong \mathbb D, \quad  w=u + i v \mapsto \R u + \R (-v).
 $$
For the isotropic matrix $\ell= \kzxz{0}{-1}{0}{0} \in L$ and $\ell' =\kzxz {0} {0} {\frac{1}N} {0} \in V$ with $(\ell, \ell')=1$. We also have the associated tube domain
$$
\mathcal H_{\ell, \ell'} =\{ \kzxz {z_1} {0} {0} {-z_2} : \,  y_1 y_2  >0\} ,\quad y_i=\hbox{Im}(z_i),
$$
together with
$$
w: \mathcal H_{\ell, \ell'} \rightarrow \mathcal L,  \quad w( \kzxz {z_1} {0} {0} {-z_2}) = \kzxz {z_1} {-2 z_1 z_2} { \frac{1}{2}} { -z_2}.
$$
This gives an isomorphism $\mathcal H_{\ell, \ell'}\cong \mathcal L/\C^\times$. We also identity $\H^2 \cup (\H^-)^2$ with $\mathcal H_{\ell, \ell'}$ by
$$\phi: z=(z_1,z_2)\mapsto\kzxz{\frac{z_1}{2}}{0}{0}{\frac{z_2}{2}}.$$ Note that we use this identification in order to have the following compatibility property and it is also the identification used in the computation of Borcherds products. The following is a special case of \cite[Proposition 3.1]{YYCM}

\begin{proposition} Define
$$
w_2: \,  \H^2 \cup (\H^-)^2  \rightarrow \mathcal L, \quad w_2(z_1, z_2) = w\circ\phi(z_1,z_2)=\kzxz {\frac{z_1}{2}} {\frac{-z_1 z_2}{2}} {\frac{1}{2}} {\frac{-z_2}{2}}.
$$
 Then the composition $pr \circ w_2$ gives an isomorphism between $\H^2 \cup (\H^-)^2$ and $\mathbb D$. Moreover, $w_2$ is $H(\R)$-equivariant, where $H(\R)$  acts on $\H^2 \cup (\H^-)^2$ via  the usual linear fraction:
$$
(g_1, g_2)(z_1, z_2) =(g_1(z_1), g_2(z_2)),
$$
and acts on $\mathcal L$ and $\mathbb D$ naturally via its action on $V$. Moreover, one has
\begin{equation} \label{eq:linebundle}
(g_1, g_2) w_2(z_1, z_2) = \frac{j(g_1, z_1) j(g_2, z_2)}{\nu(g_1, g_2)} w_2(g_1(z_1), g_2(z_2)),
\end{equation}
where $\nu(g_1, g_2) =\det g_1 =\det g_2$ is the spin character of $H=\Gspin(V)$, and
$$
j(g_1,g_2, z_1, z_2) = j(g_1, z_1) j(g_2, z_2)=(c_1 z_1+d_1)(c_2 z_2 +d_2)
$$
is the automorphy factor (of weight $(1, 1)$).
\end{proposition}

Let $\Gamma =\Gamma(2)$ in this paper, and let
$$
K(2)=\{g \in \GL_2(\hat\Z):\,  g \equiv \kzxz {*} {0} {0} {1} \pmod 2\}
$$
and $K=K_{\Gamma(2)} =(K(2) \times  K(2))\cap H(\A_f)$. Let $X_K = H(\Q) \backslash \mathbb D \times H(\A_f)/K$ be the associated Shimura variety, then one has
$$
X_K \cong  Y(2) \times Y(2),
$$
where $Y(2) =\Gamma(2) \backslash \H$. We recall that $X(2)=Y(2) \cup \{ 0, 1, \infty\}$ is the compactified modular curve. Under this identification, the tautological line bundle $\mathcal L_K= H(\Q) \backslash \mathcal L \times H(\A_f)/K$ over $X_K$ becomes the line bundle of two variable modular forms of weight $(1, 1)$ for $\Gamma(2) \times \Gamma(2)$. Notice that the dual lattice of $L$ is $L'=\frac{1}2 L$. For each $\mu \in  L'/L$ and $m \in Q(\mu)+L$, the associated special divisor $Z(m, \mu)$ is given by
$$
Z(m, \mu) =(\Gamma(2) \times \Gamma(2))\backslash \{ (z_1, z_2):\, w_2(z_1, z_2) \perp x \hbox{ for some } x \in \mu+L, Q(x)=m\}
$$
Let $S_L=\oplus_{\mu \in L'/L} \phi_\mu \subset S(V_{\A_f})$  with $\phi_\mu=\cha(\mu+\hat L)$, and let $\omega_L$ be the Weil representation  of $\SL_2(\Z)$ on $S_L$ as in \cite{YYCM}. Let $M_{0, \omega_L}^!$ be the space of $S_L$-valued weakly holomorphic modular forms of $\SL_2(\Z)$ of weight $0$ and representation $\omega_L$, i.e.,  holomorphic functions
$$
f: \H \rightarrow S_L
$$
such that  $f(\gamma \tau) =\omega_L(\gamma) f(\tau)$ and $f$ is meromorphic at the cusp $\infty$. It has Fourier expansion
\begin{equation}
f(\tau) =\sum_{\mu \in  L'/L} \sum_{\substack{m \in \Q\\ m \gg -\infty}} c(m, \mu) q^m \phi_\mu.
\end{equation}
The following  is a special case of Borcherds' far reaching lifting theorem (\cite[Theorem 13.3]{Borcherds98}, see also \cite[Theorems 2.1 and 2.2]{YYCM}).

\begin{theorem} \label{theo:BorcherdsLifting} Let $f(\tau) \in M_{0, \omega_L}^!$ be as above and assume $c(m, \mu) \in \Z$ for $m <0$. Write $z=(z_1, z_2)\in \H^2$. Then there is a meromorphic modular form $\Psi(z,f)$  of  weight $(k, k)$ for $\Gamma(2) \times \Gamma(2)$ with  $k={c(0, 0)/2}$  (with some characters) such that
\begin{enumerate}
\item  one has on the product of open modular curves $Y(2) \times Y(2)$
$$
\hbox{Div} (\Psi(z,  f)^2 )= Z(f) =\sum_{m>0, \mu\in  L'/L} c(-m, \mu) Z(m, \mu),
$$
here the multiplicity of $Z(m, \mu)$ is $2$ if $2\mu=0$ in $L'/L$ and $1$ if $2\mu\neq 0$ in $L'/L$.
\item Near each cusp $\Q\ell$ of $X_K$ ($\ell \in V$ with $(\ell, \ell)=0$), $\Psi(z, f)$ has a product expansion of the form:
$$
\Psi(z,  f)=C e((z, \rho(W, f))) \prod_{\substack{\gamma \in M' \\ (\gamma, W) >0}}
  \prod_{\substack{ \mu \in L_0'/L\\ p(\mu) \in \gamma+M}}\left[ 1 -e((\gamma, z)+(\mu, \ell'))\right]^{c(-Q(\gamma), \mu)}.
$$
Here $C$ is a constant with absolute value
\begin{equation} \label{eq:constant}
\left| \prod_{\delta \in \Z/N, \delta \ne 0}(1-e(\frac{\delta}N))^{\frac{c(0, \frac{\delta}N \ell)}2}\right|.
\end{equation}
Here the unexplained notation are the same as in  \cite[Section 2]{YYCM}, and will be defined in the special case we considered below.
\end{enumerate}
\end{theorem}

\subsection{ $\lambda(z_1) -\lambda(z_2)$ as a  Borcherds lifting}

We choose explicit representatives of cosets $L'/L=\{ \mu_j:\, 0 \le j \le 15\}$ with

$$\mu_0=\zxz{0}{0}{0}{0}, \mu_1=\zxz{0}{0}{\half}{\half}, \mu_2=\zxz{0}{0}{\half}{0},\mu_3=\zxz{0}{0}{0}{\half},$$

$$\mu_4=\zxz{\half}{0}{0}{0}, \mu_5=\zxz{\half}{0}{0}{\half}, \mu_6=\zxz{0}{\half}{0}{0}, \mu_7=\zxz{\half}{\half}{0}{0},$$

$$\mu_8=\zxz{0}{\half}{\half}{0},  \mu_9=\zxz{\half}{\half}{\half}{\half}, \mu_{10}=\zxz{0}{\half}{0}{\half}, \mu_{11}=\zxz{0}{\half}{\half}{\half}, $$

$$\mu_{12}=\zxz{\half}{0}{\half}{0}, \mu_{13}=\zxz{\half}{\half}{\half}{0}, \mu_{14}=\zxz{\half}{0}{\half}{\half}, \mu_{15}=\zxz{\half}{\half}{0}{\half}.$$

Take the cusp $\ell_\infty=-e_{12}  \in L $ and choose  $\ell_\infty'= \frac{1}2 e_{21} \in  L'$ with $(\ell_\infty, \ell_\infty') =1$. Here $e_{ij}$ is the $2\times 2$ matrix with $(ij)$ entry $1$ and other entries $0$. Let $M_\infty=(\Q \ell_\infty + \Q \ell_\infty')^\perp \cap L=\kzxz {\Z} {0} {0} {\Z}$. Take $\ell_{M_\infty}=e_{11}$ and $\ell_{M_\infty}'= \frac{1}{2} e_{22} \in M_\infty'$   with $(\ell_{M_\infty}, \ell_{M_\infty}') =1$. Furthermore, take $P_\infty=M_\infty \cap (\Q \ell_{M_\infty} + \Q \ell_{M_\infty}')^\perp =0$ in this case. Finally, take
$$
L_{0\infty}'= \{ x \in L':\,  (x, L) \subset 2\Z\},
$$
then
it is easy to check that \begin{equation}
L_{0\infty}'/L=\{\mu_0, \mu_3, \mu_4, \mu_5,
\mu_6, \mu_7, \mu_{10}, \mu_{15} \}  \label{L}.
\end{equation}
Let  $p_\infty: L_{0\infty}'/L \rightarrow  M_\infty'/M_\infty$ be  the projection defined in \cite[(2.12)]{YYCM}. Then
\begin{equation} \label{eq:p}
\begin{split}
p_\infty(\mu_0)=p_\infty(\mu_6)= 0,& \quad p_\infty(\mu_3)=p_\infty(\mu_{10}) =\frac{1}2 e_{22},\\
\quad p_\infty(\mu_4) =p_\infty(\mu_7)=\frac{1}2e_{11}, &\quad  p_\infty(\mu_5) =p_\infty(\mu_{15}) =\frac{1}2I_2.
\end{split}
\end{equation}
Here $I_2$ is the identity matrix of rank 2.

This explains most of the unexplained notations in the last theorem except the Weyl chamber $W$ and Weyl vector $\rho(W, f)$, both of which depend on the choice of $f$ (we refer to \cite[Section 2]{YYCM} for a brief review of them).  We first determine the Borcherds products of the constant weakly holomorphic vector valued modular forms which will be used later and are of independent  interest.

\begin{proposition} \label{prop:ConstantLifting} Let $M_{0, \omega_L}^{!, 0}$ be the subspace of $M_{0, \omega_L}^!$ consisting of constant vectors $f= \sum c(0, \mu_i) \phi_{\mu_i}$. Then

\begin{enumerate}
\item The space $M_{0, \omega_L}^{!, 0}$  is of dimension $5$ with a basis $\{ f_j:\, 0 \le j \le 4\}$, where
\begin{align*}
f_0&=\phi_{\mu_0} + \phi_{\mu_9} +\phi_{\mu_{10}} +\phi_{\mu_{12}} ,
\\
f_1&=\phi_{\mu_1} + \phi_{\mu_7} -\phi_{\mu_{10}} -\phi_{\mu_{12}} ,
\\
f_2&=\phi_{\mu_2} - \phi_{\mu_6} -\phi_{\mu_{7}} +\phi_{\mu_{12}} ,
\\
f_3&=\phi_{\mu_3} + \phi_{\mu_6} -\phi_{\mu_{9}} -\phi_{\mu_{12}} ,
\\
f_4&=\phi_{\mu_4} + \phi_{\mu_6} +\phi_{\mu_7}-\phi_{\mu_9}-\phi_{\mu_{10}} -\phi_{\mu_{12}} .
\end{align*}

\item Their Borcherds liftings are given by
\begin{align*}
\Psi(z, 24f_0)&=2^{12} \Delta(z_2)  \Delta(2z_1) (1- \lambda(z_1))^3=\Delta(z_2)  \omega_2(z_1)\Delta(z_1)(1- \lambda(z_1))^3,
\\
\Psi(z,24 f_1)&=\frac{\omega_2(z_2) (1-\lambda(z_2))^3}{\omega_2(z_1) (1-\lambda(z_1))^3}=\lambda(z_1)(1-\lambda(z_1))^{-2} \lambda(z_2)^{-1} (1-\lambda(z_2))^2,
\\
\Psi(z, 24f_2)&=2^{12}\frac{1}{\omega_2(z_1) \omega_2(z_2)^2 (1-\lambda(z_2))^3}=-\lambda(z_1) (1-\lambda(z_1)) \lambda(z_2)^2 (1-\lambda(z_2))^{-1},
\\
\Psi(z, 24f_3)&=-2^{-12}\omega_2(z_1)^{2} \omega_2(z_2) \lambda(z_1)^{3} =\lambda(z_1) (1-\lambda(z_1))^{-2}\lambda(z_2)^{-1} (1-\lambda(z_2))^{-1},
\\
\Psi(z, 24 f_4)&=-2^{-12}\frac{\omega_2(z_2)^3 \lambda(z_2)^3 (1-\lambda(z_2))^3}{(1-\lambda(z_1))^3} =(1-\lambda(z_1))^{-3}.
\end{align*}
Here
$$
\omega_2(z) =\frac{-16}{\lambda(z) (1-\lambda(z))}
$$
is a Hauptmodul for $\Gamma_0(2)$. Finally,  $\lambda(z_i)$ and $1-\lambda(z_i)$ can be written as products of integral powers of $\Psi(z, 8f_j)$,  $j=1,2,3,4$, up to a constant multiple.

\item If  $\Psi(z_1, z_2)$ is a two variable meromorphic modular form for $\Gamma(2) \times \Gamma(2)$ of weight $(k, k)$ with divisor solely supported on the boundary. Then
    $$
    \Psi(z_1, z_2)^3 = C\Psi(z_1, z_2, f)
    $$
is a Borcherds lifting of some $f \in M_{0,\omega_L}^{!, 0}$ up to a constant multiple.
\end{enumerate}

\end{proposition}
\begin{proof} (sketch) Write $n(b) =\kzxz {1} {b} {0} {1}$ and $w=\kzxz {0} {-1} {1} {0}$. Then  $f = \sum c_f(0, \mu_i) \phi_{\mu_i} \in M_{0, \omega_L}^{!}$ if and only if
$$
f (\tau+1)  =\omega_L(n(1)) f  \quad \hbox{ and } \quad  f(-1/\tau) =\omega_L(w)f,
$$
i.e.,
$$
\sum c_f(0, \mu_i) (e(Q(\mu_i))-1) \phi_{\mu_i} =0
$$
and
$$
\sum c_f(0, \mu_i) \omega_L(w)\phi_{\mu_i}= \sum c_f(0, \mu_i) \omega_L(w)\phi_{\mu_i}.
$$
Now a simple (tedious)  linear algebra calculation gives (1). To find the Borcherds products of $\Psi(f)$ for $f \in M_{0, \omega_L}^{!,0}$ , notice first that since $f$ has no negative terms, there is only one Weyl chamber
$$
W=\{ \R \kzxz {a} {0} {0} {-1}:\, a >0\}.
$$
So $\gamma =\kzxz {-\frac{m}2} {0} {0} {\frac{n}2} \in M'$ satisfies $(\gamma, W) >0$ if and only if $m, n\ge 0$ are integers  not both equal to zero. The Weyl vector is given  by,  via explicit calculation using the formula \cite[(2.17)]{YYCM},
$$
\rho(W, f)=\frac{1}{48} \kzxz {-c_2} {0} {0} {c_1},
$$
with
\begin{align*}
c_2&=2c_f(0, \mu_0) + 2c_f(0, \mu_6) -c_f(0, \mu_4) -c_f(0, \mu_7),
\\
c_1&=c_f(0, \mu_0) + c_f(0, \mu_6) +c_f(0, \mu_4) +c_f(0, \mu_7).
\end{align*}
Finally, the constant $C=C(f)$ can be chosen as $C(f) =2^{\frac{c_f(0, \mu_6)}{2}}$. So one has by Theorem \ref{theo:BorcherdsLifting}
\begin{align} \label{eq:LambdaConstant}
\Psi(z, f) &= 2^{\frac{c_f(0, \mu_6)}{2}} q_1^{\frac{c_1}{48}} q_2^{\frac{c_2}{48}}
\\
 &\quad \cdot \prod_{m=1}^\infty (1-q_2^m)^{c_f(0, 0)} (1+q_2^m)^{c_f(0,\mu_6)} (1-q_2^{m-\frac{1}2})^{c_f(0, \mu_4)} (1+q_2^{m-\frac{1}2})^{c_f(0, \mu_7)} \notag
 \\
 &\quad \cdot \prod_{n=1}^\infty (1-q_1^n)^{c_f(0, 0)} (1+q_1^n)^{c_f(0,\mu_6)} (1-q_1^{n-\frac{1}2})^{c_f(0, \mu_3)} (1+q_1^{n-\frac{1}2})^{c_f(0, \mu_{10})}. \notag
\end{align}
Now we have the following table

\begin{table}[h]
\begin{tabular}{|c|c|c|c|c|c|}
\hline \rule[-3mm]{0mm}{8mm}
{} & $f_0$ & $ f_1$ &$f_2$ &$f_3$ &$f_4$ \\
\hline $c_f(0, 0)$  &  $1$ &  $0$ &$0$ &$0$ &$0$\\
\hline  $c_f(0,\mu_6)$  &$0$ & $0$ &$-1$ &$1$ &$1$\\
\hline $c_f(0, \mu_4)$ &$0$ & $0$ &$0$ &$0$ &$1$\\
\hline $c_f(0, \mu_7)$ & $0$ &$ 1$ &$-1$ &$0$  &$1$\\
\hline $c_f(0, \mu_3)$ & $0$ &$ 0$ &$0$ &$1$ &$0$\\
\hline $c_f(0, \mu_{10})$ & $1$ &$-1$ &$0$ &$0$ &$-1$\\
\hline  $c_1$  &$1$ &$1$ &$-2$ &$1$ &$3$\\
\hline $c_2$ &$2$ &$-1$ &$-1$ &$2$ &$0$\\
\hline $C(f)$ &$1$ &$1$ &$2^{-\half}$ &$2^{\half}$ &$2^{\half}$\\
\hline
\end{tabular}
\end{table}

Recall the product formula (\ref{eq:lambda}), and
$$
1-\lambda(\tau)=\frac{1}{16}q_\tau^{-\frac{1}{2}}\prod_{n\geq1}(\frac{1+q_\tau^{n-\frac{1}{2}}}{1+q_\tau^n})^8,
$$
we  obtain the formulas in  (2).  We rewrite the formulas as
$$
\begin{pmatrix}
\log|\Psi(z, 24f_1)| \\
\log|\Psi(z, 24f_2)|  \\
\log|\Psi(z, 24f_3)| \\
\log|\Psi(z, 24f_4)|
\end{pmatrix}
= \begin{pmatrix}
1 &-2 &-1 &2\\
1 &1 &2 &-1\\
1 &-2 &-1 &-1 \\
0 &-3 &0 &0
\end{pmatrix}
\begin{pmatrix}
\log|\lambda(z_1)| \\
\log|1-\lambda(z_1)|  \\
\log|\lambda(z_2) | \\
\log|1-\lambda(z_2)|
\end{pmatrix}.
$$
The matrix is invertible with  inverse
$$\frac{1}{3}\begin{pmatrix}
1 &1 &1 &-1\\
0 &0 &0 &-1\\
0 &1 &-1 &1\\
1 &0 &-1 &0
\end{pmatrix}.$$
This proves the last claim of (2): $\lambda(z_1)$ and $1-\lambda(z_1)$ are also Borcherds products.

To prove (3), notice that $\Psi_1(z)=\Psi(z)/\Psi(z, 2k f_0)$ is a two variable modular function with divisor supported on the boundary of $X(2) \times X(2)$. It is easy to check that  such a function is of the form
$$
\Psi_1(z)=C_1 \lambda(z_1)^{m_1} (1-\lambda(z_2))^{m_2} \lambda(z_2)^{m_3} (1-\lambda(z_2))^{m_4}
$$
with $m_i \in \Z$. Indeed, fix $z_2 \in  \H$, then $\Psi_1(z_1, z_2)$ is a meromorphic function on $X(2)$ with poles and zeros only at the cusps $\{0, 1, \infty\}$, so
$$
\Psi_1(z_1, z_2)=C(z_2) \lambda(z_1)^{m_1} (1-\lambda(z_1))^{m_2},
$$
where $C(z_2)$ is a meromorphic function of $z_2$ which can only have zeros or poles at the cusps of $X(2)$. So
$$
C(z_2)=C_1 \lambda(z_2)^{m_3} (1-\lambda(z_2))^{m_4}.
$$
So  there are some integers $b_i$   by (2)  such that
$$
\Psi_1(z)^3 =C \Psi(z, f_1)^{b_1} \Psi(z, f_2)^{b_2} \Psi(z, f_3)^{b_3} \Psi(z, f_4)^{b_4}.
$$
This proves the proposition.

\end{proof}

\begin{corollary}\label{lambda}We have the following Borcherds lifting:
\begin{eqnarray*}
-\lambda(z_1)&=&\Psi(z, 8(f_1+f_2+f_3-f_4))=\Psi(z,8(\phi_{\mu_1}+\phi_{\mu_2}+\phi_{\mu_3}-\phi_{\mu_4}-\phi_{\mu_6}-\phi_{\mu_7})),\\
1-\lambda(z_1)&=&\Psi(z,-8f_4)=\Psi(z,-8(\phi_{\mu_4} + \phi_{\mu_6} +\phi_{\mu_7}-\phi_{\mu_9}-\phi_{\mu_{10}} -\phi_{\mu_{12}})),\\
-\lambda(z_2)&=&\Psi(z, 8(f_2-f_3+f_4))=\Psi(z,8(\phi_{\mu_2}-\phi_{\mu_3}+\phi_{\mu_4}-\phi_{\mu_6}-\phi_{\mu_{10}}+\phi_{\mu_{12}})),\\
1-\lambda(z_2)&=&\Psi(z,8(f_1-f_3))=\Psi(z,8(\phi_{\mu_1}-\phi_{\mu_3}-\phi_{\mu_6}+\phi_{\mu_7}+\phi_{\mu_9}-\phi_{\mu_{10}})).
\end{eqnarray*}
\end{corollary}

\begin{proposition} \label{prop6.2} Let
\begin{align*}
f(\tau) &=-8 \sum_{\gamma \in \Gamma(2)\backslash \SL_2(\Z)} \lambda|\gamma \omega_L(\gamma)^{-1} \phi_{\mu_5}
        -4 f_0 + 4 f_1+ 4 f_2 -4 f_3 -4 f_4
   \\
    &= q^{-\frac{1}2} \phi_{\mu_5} + \sum_{m \ge 0, \mu \in L'/L} c_f(m, \mu) q^m \phi_\mu  \in M_{0, \omega_L}^!
\end{align*}
be a weakly holomorphic vector valued modular form induced from $\lambda$ with slightly modification as indicated. Then the principal part of $f$ is $q^{-\frac{1}2}\phi_{\mu_5}$ and  (with a suitable choice of the constant)
$$
\Psi(z, f) =\lambda(z_1) -\lambda(z_2),
$$
and
\begin{align*}
c_f(0, \mu_3)&=c_f(0, \mu_4)=c_f(0, \mu_6)=c_f(0, \mu_7)=c_f(0, \mu_{10})=-8,
\\
 c_f(0,  \mu_2) &=c_f(0,\mu_{9}) =8,
 \\
 c_f(0, \mu_j) &=0 \hbox{  for all other $j$}.
 \end{align*}
\end{proposition}
\begin{proof} Let
$$
f_L(\tau)=-8 \sum_{\gamma \in \Gamma(2)\backslash \SL_2(\Z)} \lambda|\gamma \omega_L(\gamma)^{-1} \phi_{\mu_5}
$$
be the weakly holomorphic vector valued modular form induced from $\lambda$. Then one has the following Fourier expansion
$$
f_L(\tau)= q^{-\half}\phi_{\mu_5}+ \sum_{j=0}^{15}\sum_{r \in \Z_{\ge 0} }c_L(\frac{r}2,\mu_j)q^{\frac{r}2}   \phi_{\mu_j},
$$
with
\begin{align*}
c_L(0, \mu_0)&=c_L(0, \mu_2)=c_L(0, \mu_6)=c_L(0, \mu_9)=4,
\\
 c_L(0,  \mu_1) &=c_L(0,\mu_{3})=c_L(0,  \mu_4) =c_L(0,  \mu_7) =c_L(0,  \mu_{10}) =c_L(0,  \mu_{12}) =-4.
 \end{align*}
So $Z(f_L) =Z(\frac{1}2, \mu_5) =Y_{\Gamma(2)}^\Delta$ is the diagonal cycle in the open Shimura variety $X_K=Y(2) \times Y(2)$.


We first look at the Borcherds lifting of $f_L$ around the cusp $\Q\ell_\infty$ with $\ell_\infty =-e_{12}$ as at the beginning of this section, and $\ell_\infty'=\frac{1}2 e_{21}$.
The Grassmannian $\Gr(M)$ is cut into two Weyl chambers by $Z(f_L)$, the one whose closure contains $\ell_{M_\infty}=e_{11}$ is
$$
W=\{ \R \kzxz {a} {0} {0} {-1}:\, a >1\}.
$$
The associated Weyl vector is $\rho(W, f_L)=0$.  For a vector $\gamma=\kzxz {-\frac{m}2} {0} {0} {\frac{n}2}  \in M'$, $(\gamma, W) >0$ if and only if
$$
m+n \ge 0, \quad n >0,\quad  m \hbox{ and } n \hbox{ are not both } 0.
$$
Finally, the identification between $\H^2\cup (\H^-)^2$ and $\mathcal H_{\ell, \ell'}$ is  $\phi(z_1, z_2) = \kzxz {\frac{z_1}{2}} {0} {0} {-\frac{z_2}2}$.  So Theorem \ref{theo:BorcherdsLifting}  gives in this special case
\begin{equation}\label{Borcherds1}
\Psi(z_1, z_2, f_L)=q_2^{\frac{1}2} (1- q_1^{\frac{1}2} q_2^{-\frac{1}2})
 \prod_{m, n \ge 1} \prod_{\substack{ \mu \in L_{0\infty}'/L \\ p_\infty(\mu) =\kzxz {-\frac{m}2} {0} {0} {\frac{n}2} +M_\infty}}(1-q_1^{\frac{n}2} q_2^{\frac{m}2} e((\mu, \ell_\infty')))^{c(\frac{mn}2, \mu)}.
\end{equation}

Let $\Psi(S(z_1), S(z_2), f_L)$ be the Borcherds product expansion of $\Psi(z, f_L)$ at the
cusp $\ell_0= w \ell_\infty w^{-1} = e_{21}$. So $\ell_0'=w \ell_\infty' w^{-1} =- \frac{1}2 e_{12}$. Theorem \ref{theo:BorcherdsLifting}, together with the same calculation as above, gives the Borcherds product expansion:
\begin{equation}\label{Borcherds2}
\Psi(S(z_1), S(z_2), f_L)=q_2^{\frac{1}2}(1- q_1^{\frac{1}2} q_2^{-\frac{1}2})
 \prod_{m, n \ge 1} \prod_{\substack{ \mu \in L_{00}'/L \\ p_0(\mu) =\kzxz {\frac{m}2} {0} {0} {-\frac{n}2} +M_0}}
   (1-q_1^{\frac{m}2} q_2^{\frac{n}2} e((\mu, \ell_0'))^{c(\frac{mn}2, \mu)},
\end{equation}
 where
$M_0=(\Q \ell_0 + \Q \ell_0')^\perp \cap L=\kzxz {\Z} {0} {0} {\Z}
$, and

$$
L_{00}'=\{ x \in L': \,  (x, \ell_0') \in 2\Z\}=\{\mu_0, \mu_1, \mu_2, \mu_3,
\mu_4, \mu_5, \mu_{12}, \mu_{14} \} +L,
$$
and
$$
p_0:  L_{00}'/L \rightarrow M_0'/M_0,
$$
is the associated projection similar to (\ref{eq:p}).

Now we assume $$f'=\sum_{i=0}^{15}a_i\phi_{\mu_i}$$ is a constant vector valued modular form such that $\Psi(z_1,z_2,f_L+f')=\lambda(z_1)-\lambda(z_2)$.  From the expression (\ref{Borcherds1}) and (\ref{Borcherds2}) and comparing the leading terms, we want the Borcherds lifting of $f'$ at the cusp $\ell_{\infty}$ to have  Weyl vector $\zxz{0}{0}{0}{-\frac{1}{2}}$ and the Borcherds lifting of $f'$ at the cusp $\ell_0$ to have Weyl vector $\zxz{0}{0}{0}{-\frac{1}{2}}$. Note that

$$f'_{M_{\infty}}=(a_0+a_6)\phi_{\mu_0}+(a_4+a_7)\phi_{\mu_4},$$
$$f'_{M_{0}}=(a_0+a_2)\phi_{\mu_0}+(a_1+a_3)\phi_{\mu_3},$$
and
$$\ell_{M_{\infty}}=e_{11},\  \ell_{M_{\infty}}'=\half e_{22}, \ \ell_{M_{0}}=e_{22},\ \ell_{M_{0}}'=\half e_{11}.$$
Then  by the Weyl vector formula in \cite[(2.15)-(2.18)]{YYCM},  we have

$$\rho_{\ell_{M_{\infty}}'}=\frac{a_0+a_4+a_6+a_7}{24}=-1,$$

$$ \rho_{\ell_{M_{\infty}}}=-\frac{a_0+a_6}{24}+\frac{a_4+a_7}{48}=\frac{1}{2},$$

$$\rho_{\ell_{M_{0}}'}=\frac{a_0+a_1+a_2+a_3}{24}=0,$$

$$\rho_{\ell_{M_{0}}}=-\frac{a_0+a_2}{24}+\frac{a_1+a_3}{48}=0$$
Since $c(0,f_L)=4$, we also require $a_0=-4$ in order to get a modular function. Note that $T\circ f'=f'$ and $S\circ f'=f'$ is the basic requirement for $f'$ to be a vector valued modular form. These linear equations have a unique solution which can be written as
$$f'=-4 f_0 + 4 f_1+ 4 f_2 -4 f_3 -4 f_4.$$
Now it is easy to see that
\[\Psi(z_1,z_2,f_L+f')/(\lambda(z_1)-\lambda(z_2))\]
is holomorphic  on $X(2)\times X(2)$ and is thus constant. By comparing the coefficients, we get that  $\Psi(z_1,z_2,f_L+f')=\lambda(z_1)-\lambda(z_2)$.
\end{proof}

For the convenience of the reader interested in calculation, we record an  explicit expansion for $f(\tau)$ as the following proposition.

\begin{proposition}  \label{prop:Explicitf} Let
$$
f(\tau) = \sum_{j=0}^{15} f(\tau, \mu_j) e_{\mu_j}
$$
be the weakly modular form in Proposition \ref{prop6.2}. Then
\begin{eqnarray*}
f(\tau,\mu_0)&=&1024q+24576q^2+307200q^3+2686976q^4+18561024q^5+\cdots\\
f(\tau,\mu_1)&=&-1024q-24576q^2-307200q^3-2686976q^4-18561024q^5+\cdots\\
f(\tau,\mu_2)&=&8+1024q+24576q^2+307200q^3+2686976q^4+18561024q^5+\cdots\\
f(\tau,\mu_3)&=&-8-1024q-24576q^2-307200q^3-2686976q^4-18561024q^5+\cdots\\
f(\tau,\mu_4)&=&-8-1024q-24576q^2-307200q^3-2686976q^4-18561024q^5+\cdots\\
f(\tau,\mu_5)&=&q^{-\frac{1}{2}}+148q+5570q^2+92120q^3+940415q^4+7239908q^5+\cdots\\
f(\tau,\mu_6)&=&-8+1024q+24576q^2+307200q^3+2686976q^4+18561024q^5+\cdots\\
f(\tau,\mu_7)&=&-8-1024q-24576q^2-307200q^3-2686976q^4-18561024q^5+\cdots\\
f(\tau,\mu_8)&=&128q^{\frac{1}{2}}+5632q^{\frac{3}{2}}+91904q^{\frac{5}{2}}+941056q^{\frac{7}{2}}+7238272q^{\frac{9}{2}}+\cdots\\
f(\tau,\mu_9)&=&8+1024q+24576q^2+307200q^3+2686976q^4+18561024q^5+\cdots\\
f(\tau,\mu_{10})&=&-8-1024q-24576q^2-307200q^3-2686976q^4-18561024q^5+\cdots\\
f(\tau,\mu_{11})&=&-128q^{\frac{1}{2}}-5632q^{\frac{3}{2}}-91904q^{\frac{5}{2}}-941056q^{\frac{7}{2}}-7238272q^{\frac{9}{2}}+\cdots\\
f(\tau,\mu_{12})&=&-1024q-24576q^2-307200q^3-2686976q^4-18561024q^5+\cdots\\
f(\tau,\mu_{13})&=&-128q^{\frac{1}{2}}-5632q^{\frac{3}{2}}-91904q^{\frac{5}{2}}-941056q^{\frac{7}{2}}-7238272q^{\frac{9}{2}}+\cdots\\
f(\tau,\mu_{14})&=&128q^{\frac{1}{2}}+5632q^{\frac{3}{2}}+91904q^{\frac{5}{2}}+941056q^{\frac{7}{2}}+7238272q^{\frac{9}{2}}+\cdots\\
f(\tau,\mu_{15})&=&128q^{\frac{1}{2}}+5632q^{\frac{3}{2}}+91904q^{\frac{5}{2}}+941056q^{\frac{7}{2}}+7238272q^{\frac{9}{2}}+\cdots
\end{eqnarray*}

\end{proposition}

\section{The Big CM value formula and proofs of the main results} \label{sect:BigCM}

\subsection{Big CM cycles}
Let $\kay_i=\Q(\sqrt{d_i})$, $i=1, 2$, be two imaginary quadratic fields of fundamental discriminants $d_i$, and $E=\kay_1\kay_2 =\Q(\sqrt{d_1}, \sqrt{d_2})$. Let $F=\Q(\sqrt D)$ be the real quadratic field with $D=d_1 d_2$. We assume in this paper $(d_1, d_2)=1$. Let $W=E$ with $F$-quadratic form $Q_F(x) =\frac{ 2 x \bar x}{\sqrt D}$, and let $W_\Q=W$ with $\Q$-quadratic form $Q(x) =\tr_{F/\Q}(Q_F(x))$.  Let $\sigma_1=1$ and $\sigma_2 =\sigma$ be two real embeddings of $F$ with $\sigma_i(\sqrt D) =(-1)^{i-1}\sqrt D$. Then
 $W$ has signature $(0, 2)$ at $\sigma_2$ and $(2, 0)$  at $\sigma_1$ respectively, and so $W_\Q$ has signature $(2,2)$.
 Choose a $\Z$-basis of $\OO=\OO_E$ as follows
 $$
 e_1=1\otimes 1, \quad e_2 = \frac{-d_1+\sqrt{d_1}}2 =\frac{-d_1+\sqrt{d_1}}2 \otimes 1, \quad  e_3= \frac{d_2+\sqrt{d_2}}2 =1\otimes \frac{d_2+\sqrt{d_2}}2, \quad e_4= e_2 e_3.
 $$
We will drop $\otimes$ when there is no confusion.  Then it is easy to check that
\begin{equation} \label{eq:SpaceIdentification}
(W_\Q, Q_\Q) \cong (V, Q)=(M_2(\Q),  2\det), \sum x_i e_i \mapsto \kzxz {x_3} {x_1} {x_4} {x_2}.
\end{equation}
We will identify $(W_\Q, Q_\Q)$ with the quadratic space $(V, Q) =(M_2(\Q), 2\det)$. Under this identification, $L=M_2(\Z)$ becomes $\OO_E$.

Let $T$ be the torus over $\Q$ given by  (see \cite{HYbook}) and also \cite[Section 6]{BKY12})
$$
T(R) =\{ (t_1, t_2)\in (\kay_1\otimes_\Q R)^\times \times (\kay_2\otimes_\Q R)^\times:\,  t_1 \bar{t}_1 = t_2 \bar{t}_2 \},
$$
for any $\Q$-algebra $R$. Then there is a    map from $T$ to $\SO(W)$ given by  $(t_1, t_2) \mapsto {t_1}/{\bar{t}_2}$.
Define the embedding
\begin{equation} \label{eq:iota}
\iota_j: \kk_j \rightarrow  M_2(\Q), \quad  \iota_j(r) ( e_{j+1}, e_1)^t = ( r e_{j+1}, r e_1)^t.
\end{equation}
Then $\iota= (\iota_1, \iota_2)$ gives the embedding from $T$ to $H$, and one has the commutative diagram:
$$
\begin{matrix}
1 &\rightarrow &\G_m &\rightarrow  &T &\rightarrow  &\Res_{F/\Q} \SO(W) &\rightarrow &1
\\
 &             &\downarrow &       &\downarrow  &  &\downarrow &  &
 \\
1 &\rightarrow &\G_m  &\rightarrow  &H &\rightarrow  & \SO(V)&\rightarrow &1
\end{matrix}
$$
This shows that $T$ is a maximal torus of $H$. Since $W_{\sigma_2}$ is of signature $(0, 2)$, it gives two points $z_{\sigma_2}^\pm$ in $X_K$, a big CM point in the sense of \cite{BKY12}. The associated big CM cycles are
$$
Z(W, z_{\sigma_2})=\{ z_{\sigma_2}\} \times (T(\Q) \backslash T(\A_f)/K_T)
$$
where  $z_{\sigma_2} =z_{\sigma_2}^+$,
$$K_T=\{(t_1, t_2) \in T(\A_f):\,  (\iota_1(t_1), \iota_2(t_2)) \in K(\Gamma(2)) \times K(\Gamma(2))\}.$$
This cycle is defined over $\sigma_2(E)=E$. Let $Z(W)$ be the formal sum of all its Galois conjugates (four of them), which is a big CM cycle of  $X_K$ defined over $\Q$. We refer to \cite{BKY12} for details.

\begin{lemma} \label{lem:ClassIsomorphism} Let $K_j=1+2 \hat{\OO}_{{\kay_j}} \subset \hat{\OO}_{{\kay_j}}^\times$, then $\Cl(K_j) = {\kay_j}^\times \backslash E_{j,f}^\times/K_j$ is the ray class group $\Cl({\kay_j}, 2)$ of $E$ of modulus $2\OO_{{\kay_j}}$.
Moreover, one has an isomorphism
$$
p': T(\Q) \backslash  T(\A_f) /K_T \rightarrow  \Cl(\kay_1, 2) \times \Cl(\kay_2, 2), \quad p'([t_1, t_2]) =([t_1], [t_2]).
$$
\end{lemma}
\begin{proof}First notice that $K_j=\iota_j^{-1} (K(2))$.  Let $H(\kay_i, N)$ be the ray class field of $\kay_i$ of modulus $N\OO_{\kay_i}$ for $i=1, 2$ and $N=1, 2$. By the same argument as that in \cite[Lemma 3.8]{YYCM}, it is enough to show $H(\kay_1, 2) \cap H(\kay_2, 2)=\Q$. Since $(d_1, d_2) =1$, we may assume $2\nmid d_2$.  In particular, $2$ is unramified in  $H(\kay_2, 1)$. Every rational prime $p\ne 2$ is unramified in $H(\kay_1, 2) \cap H(\kay_2, 2)$ and thus in $H(\kay_1, 2) \cap H(\kay_2, 1)$. Therefore $H(\kay_1, 2)\cap H(\kay_2, 1)$ is unramified over $\Q$ at every finite prime and  is thus equal to $\Q$. Now we have
$$
[H(\kay_1, 2) \cap H(\kay_2, 2):\Q] =[H(\kay_1, 2) \cap H(\kay_2, 2):H(\kay_1, 2) \cap H(\kay_2, 1)],
$$
which divides  $[H(\kay_2, 2): H(\kay_2, 1)]$. Since $[H(\kay_2, 2): H(\kay_2, 1)]=1$ or  $3$ by Lemma \ref{lem:ClassGroup},  $H(\kay_1, 2) \cap H(\kay_2, 2)=\Q$ or $H(\kay_1, 2) \cap H(\kay_2, 2)$ is abelian over $\Q$ of degree $3$. Since  $H(\kay_1, 2)\cap H(\kay_2, 2)$ is also unramified over $\Q$ at finite primes away from $2$, we have $H(\kay_1, 2)\cap H(\kay_2, 2) \subset \Q(\zeta_{2^{n}})$, which has even degree over $\Q$, a contradiction. Therefore $H(\kay_1, 2) \cap H(\kay_2, 2)=\Q$. This proves the lemma.

\end{proof}

\subsection{Incoherent Eisenstein series and the big CM value formula of Bruinier, Kudla, and Yang}  In  this subsection, we review  the big CM value formula of Bruinier, Kudla, and Yang \cite{BKY12}.

Associated to the $F$-quadratic space  $W$ and the additive adelic character $\psi_F =\psi \circ \tr_{F/\Q}$ is a Weil representation $\omega=\omega_{\psi_F}$ of $\SL_2(\A_F)$  on $S(W(\A_F)) =S(V(\A_\Q))$. Let $\chi=\chi_{E/F} $ be the quadratic Hecke character of $F$ associated to $E/F$, then $\chi=\chi_W$ is also the quadratic Hecke character of $F$ associated to $W$, and there is a $\SL_2(\A_F)$-equivariant map
\begin{equation}
\lambda=\prod \lambda_v : S(W(\A_F)) \rightarrow  I(0, \chi), \quad  \lambda(\phi) (g) = \omega(g) \phi(0).
\end{equation}
Here $I(s, \chi) =\Ind_{B_{\A_F}}^{\SL_2(\A_F) }\chi |\,|^s$ is the principal series, whose sections (elements) are smooth functions $\Phi$  on $\SL_2(\A_F)$ satisfying the condition
$$
\Phi(n(b) m(a) g,s )= \chi(a)|a|^{s+1}\Phi(g, s), \quad b \in \A_F,  \hbox{ and  } a \in \A_F^\times.
$$
Here $B =NM$ is the standard Borel subgroup of $\SL_2$. Such a section is called factorizable if $\Phi=\otimes \Phi_v$ with $\Phi_v\in I(s, \chi_v)$. It is called standard if $\Phi|_{\SL_2(\hat{\OO}_F) \SO_2(\R)^{2}} $ is independent of $s$. For a standard  section $\Phi \in I(s, \chi)$, its associated Eisenstein series is defined  for $\Re(s) \gg 0$ as
$$
E(g, s, \Phi) = \sum_{\gamma \in B(F) \backslash \SL_2(F)} \Phi(\gamma g, s).
$$

For $\phi \in S(V_f) =S(W_f)$, let $\Phi_f $ be the standard section associated to $\lambda_f(\phi) \in I(0, \chi_f)$, here the subscript $f$ means the finite part of the adele. For each real embedding $\sigma_i: F \hookrightarrow \R$,  let  $\Phi_{\sigma_i} \in I(s, \chi_{\C/\R})=I(s, \chi_{E_{\sigma_i}/F_{\sigma_i}})$ be the unique `weight one' eigenvector of $\SL_2(\R)$ given by
$$
\Phi_{\sigma_i}(n(b)m(a) k_\theta) = \chi_{\C/\R}(a) |a|^{s+1} e^{i  \theta},
$$
for $b \in \R$, $a\in \R^\times$, and $k_\theta =\kzxz {\cos\theta} {\sin \theta} {-\sin \theta} {\cos \theta} \in \SO_2(\R)$. We define  for $\vec\tau =(\tau_1,  \tau_{2}) \in \H^{2}$
$$
E(\vec\tau, s, \phi) =  \norm(\vec v)^{-\frac12} E(g_{\vec\tau}, s, \Phi_f \otimes  (\otimes_{1 \le i \le 2}\Phi_{\sigma_i} )),
$$
where $\vec v =\hbox{Im}(\vec\tau)$, $\norm(\vec v) =\prod_i v_i$, and $g_{\vec\tau} = (n(u_i) m(\sqrt{v_i}))_{1\le i \le 2}$. It is a (non-holomorphic) Hilbert modular form of scalar weight $1$ for some congruence subgroup of $\SL_2(\OO_F)$. Following \cite{BKY12}, we further normalize
$$
E^*(\vec\tau, s, \phi) = \Lambda(s+1, \chi) E(\vec\tau, s, \phi),
$$
where $\partial_F$ is the different of $F$, $d_{E/F}$ is the relative discriminant of $E/F$, and
\begin{equation} \label{eq:L-series}
\Lambda(s, \chi) =A^{\frac{s}2} (\pi^{-\frac{s+1}2}
\Gamma(\frac{s+1}2))^{2} L(s, \chi)=\Lambda(s, \chi_{\kay_1/\Q}) \Lambda(s, \chi_{\kay_2/\Q}),  \quad A=\norm_{F/\mathbb Q}
(\partial_F d_{E/F})=D.
\end{equation}
$\Lambda(s,\chi)$ is a holomorphic function of $s$ with functional equation
\[\Lambda(s,\chi)=\Lambda(1-s,\chi),\]
and
\[\Lambda(1,\chi)=\Lambda(0,\chi)=L(0,\chi)=2^{2-\delta}\frac{h(E)}{w(E)h(F)}=4\frac{h(\kay_1)h(\kay_2)}{w(\kay_1)w(\kay_2)}\in \Q^\times,\]
where $2^\delta=[\mathcal O_E^\times:\mu(E)\mathcal O_F^\times]=1,2.$

The Eisenstein series is incoherent in the sense of Kudla (\cite{KuAnnals}). This forces $E^*(\vec\tau, 0, \phi)=0$ automatically.


\begin{proposition}(\cite[Proposition 4.6]{BKY12}) \label{prop4.3} Let $\ph \in S(V_f) =S(W_f)$.
For  a totally positive element $t\in F^\times_+$,
let $a(t,\ph)$ be the $t$-th Fourier coefficient of $E^{*,
\prime}(\vec\tau, 0, \ph)$ and write
the constant term  of  $E^{*, \prime}(\vec\tau, 0, \ph)$ as
$$\ph(0)\Lambda(0, \chi)  \log \norm(\vec v) +a_0(\ph).$$
Let
$$
\mathcal E(\tau, \ph) =  a_0(\ph) + \sum_{n \in \mathbb
Q_{>0}} a_n(\ph)\, q^n,
$$
where (for $n >0$)
$$
a_n(\ph) =\sum_{t \in F^\times_+, \, \tr_{F/\mathbb Q} t =n} a(t,\ph).
$$
Here $F_+^\times$ consists of all totally positive elements in $F$.
Then, writing $\tau^\Delta=(\tau,  \tau)\in \H^2$,
$$
E^{*, \prime}(\tau^\Delta, 0, \ph) -\mathcal E(\tau, \ph)-2\ph(0) \Lambda(0, \chi)\,
\log v
$$
is of exponential decay as  $v$ goes to infinity.
Moreover, for $n >0$
$$
a_n(\ph)  =\sum_p a_{n, p}(\ph) \log p
$$
with $a_{n, p}(\ph) \in \mathbb Q(\ph)$, the subfield of $\mathbb C$
generated by the values $\ph(x)$, $x\in V(\A_f)$.
\end{proposition}

\begin{remark} \label{rem:Constant}  There is a minor mistake in \cite[Proposition 4.6]{BKY12}) about the constant. The corrected form is
$$
E_0^{*, \prime}(\vec\tau, 0, \ph)=\ph(0)\Lambda(0, \chi)  \log \norm(\vec v) +a_0(\ph)
$$
(i.e., $a_0(\phi)$ might not be a multiple of $\phi(0)$). Indeed,  direct calculation gives
$$
E_0^*(\vec\tau, s, \ph) = \ph(0) \Lambda(s+1, \chi) (\norm(\vec v))^{\frac{s}2} -(\norm(\vec v))^{-\frac{s}2} \Lambda(s, \chi) \tilde W_{0, f}(s, \ph),
$$
where     (when $\phi$ is factorizable)
$$
\tilde W_{0, f}(s, \ph)= \prod_{\mathfrak p \nmid \infty} \tilde W_{0, \mathfrak p}(s, \phi_\mathfrak p)
= \prod_{\mathfrak p \nmid \infty} \frac{|A|_\mathfrak p^{-\frac{1}2} L_\mathfrak p(s+1, \chi)}{\gamma(W_\mathfrak p) L_\mathfrak p(s, \chi)}  W_{0, \mathfrak p}(s, \phi_\mathfrak p)
$$
is the product of re-normalized local  Whittaker functions. With this notation, one has
\begin{equation} \label{eqa0}
a_0(\phi) = - \Lambda(0, \chi)\tilde{W}_{0, f}'(0, \phi) - 2\phi(0)  \Lambda'(0, \chi).
\end{equation}
Here we used the fact that $\tilde W_{0, f}(0,\phi) =\phi(0)$, which follows from $E_0^*(\tau, 0, \phi)=0$.
\end{remark}

Notice that $a(t, \phi_\mu) =0$ automatically unless $\mu+\hat L$ represents $t$, i.e., $t -Q_F(\mu) \in  \partial_F^{-1}\OO_F$.
The following is  now a special case of    the big CM value formula of Bruinier, Kudla, and Yang (\cite[Theorem 5.2]{BKY12}).

\begin{theorem} \label{theo:BigCM}Assume $(d_1,d_2)=1$. Let $\tau_i =\frac{d_i+\sqrt{d_i}}2$, $i=1, 2$, and let
$$
f(\tau) =\sum_{\mu \in L'/L} f_\mu(\tau) \phi_\mu =\sum c(m, \mu) q^m \phi_\mu \in M_{0, \omega_L}^!
$$
be a weakly holomorphic modular form of $\SL_2(\Z)$ of  weight $0$ with respect to the Weil representation $\omega_L$
and $c(0, 0)=0$. Let $\Psi(z, f)$ be its Borcherds  lifting. Then
\begin{align*}
&-\sum_{[\mathfrak a_i] \in \Cl(\kay_i, 2)}
\big( \log |\Psi(\tau_1^{\sigma_{\mathfrak a_1}}, \tau_2^{\sigma_{\mathfrak a_2}}, f)|^4 +\log |\Psi((-\bar{\tau}_1)^{\sigma_{\mathfrak a_1}}, \tau_2^{\sigma_{\mathfrak a_2}}, f)|^4
\\
&\qquad  + \log |\Psi(\tau_1^{\sigma_{\mathfrak a_1}}, (-\bar{\tau}_2)^{\sigma_{\mathfrak a_2}}, f)|^4+\log |\Psi((-\bar{\tau}_1)^{\sigma_{\mathfrak a_1}}, (-\bar{\tau_2})^{\sigma_{\mathfrak a_2}}, f)|^4
 \big)
\\
&=2 r(d_1) r(d_2)\left(  \sum_{\substack { 0 \ne \mu \in  L'/L\\ m\ge 0 \\ m \equiv Q(\mu) \pmod 1 }} c(-m, \mu) a_m(\phi_\mu)    \right).
\end{align*}
\end{theorem}
\begin{proof} One has by \cite[Theorem 5.2]{BKY12}
$$
-\log| \Psi(Z(W), f) |^4 =C(W, K)\left(  \sum_{\substack { 0 \ne \mu \in  L'/L,\ m \ge 0 \\ m \equiv Q(\mu) \pmod 1 }} c(-m, \mu) a_m(\phi_\mu)    \right).
$$
 Here
$$
C(W, K)=\frac{ \deg(Z(W,z_{\sigma_2}^\pm))}{  \Lambda(0, \chi)}.
$$
On the other hand, one has $z_{\sigma_2} = (\tau_1, \tau_2)$ and
$$
[z_{\sigma_2}, t_1^{-1}, t_2^{-1}] =(\tau_1^{\sigma_{t_1}}, \tau_2^{\sigma_{t_2}})
$$
 by \cite[Lemma 3.4, Proposition 3.6]{YYCM}. So $-\log| \Psi(Z(W), f) |^4$ is equal to the left hand side of the identity in the theorem by Lemma \ref{lem:ClassIsomorphism}. On the other hand, Lemma \ref{lem:ClassIsomorphism} also implies
 \begin{align*}
 C(W, K) &= \frac{2 h(\kay_1, 2)h(\kay_2, 2)}{\Lambda(0, \chi_{\kay_1/\Q}) \Lambda(0, \chi_{\kay_2/\Q})}
 \\
 &= 2 \prod_{i=1}^2 \frac{2 h(\kay_i, 2)}{w_i h(\kay_i)} =2r(d_1)r(d_2).
 \end{align*}
This proves the theorem.
\end{proof}

\subsection{The norm of $\lambda(\tau_1) -\lambda(\tau_2)$}

Combining Theorem  \ref{theo:BigCM} and Lemma \ref{lem:ClassGroup}, we have

\begin{proposition} \label{prop:norm2}Assume $(d_1,d_2)=1$. Let $\tau_j=\frac{d_j +\sqrt{d_j}}2$ with $j=1, 2$, and $\delta(D) =1$ or $2$ depending on whether $D$ is odd or even.
Then
\[
\log|\norm(\lambda(\tau_1) -\lambda(\tau_2))|
= -\frac{r(d_1)r(d_2)}{2\delta(D)}
 \sum_{ \substack{ t =\frac{m+\sqrt D}4  \in \frac{1}2 \OO_F \\ |m| < \sqrt D\\t -2\mu_5\bar{\mu}_5\in\OO_F}} a(\frac{t}{\sqrt D}, \phi_{\mu_5})-\frac{r(d_1)r(d_2)}{2\delta(D)}a_0,
\]
where
\[a_0=
8a_0(\phi_{\mu_2}) -8a_0(\phi_{\mu_3})-8a_0(\phi_{\mu_4})-8a_0(\phi_{\mu_6})-8a_0(\phi_{\mu_7}) +8a_0(\phi_{\mu_{9}})-8a_0(\phi_{\mu_{10}}).\]
\end{proposition}

\begin{remark}
The appearance of $\delta(D)$ is due to the fact that there is one Galois orbit in $Z(W)$ when $D$ is odd and two when
$D$ is even.
\end{remark}

\subsection{ Computing $a(\frac{t}{\sqrt D}, \phi_{\mu_5})$}

To prove the main theorem of this paper, we need to compute $a(\frac{t}{\sqrt D}, \phi_{\mu_5})$ for $t = \frac{m+\sqrt D}4 \in \frac{1}2 \OO_F$ with  $|m|< \sqrt D$, which occupies this subsection, and $a_0(\phi_\mu)$ for various $\mu=\mu_j$, which will occupy the next subsection. We first set up some notation.  Let $\gamma(W_v) =\gamma(W_v, \ph_v)$ be the local Weil index, which is an 8-th root of unity. Then it is known (see for example \cite[(3.4)]{KuSplit})
$$
\gamma(W_{\sigma_2}) = i = \gamma(W_{\sigma_1})^{-1},  \quad \prod_{\mathfrak p\nmid \infty} \gamma(W_{\mathfrak p})=1.
$$
We normalize
$$
W_{t, \mathfrak p}^*(s, \phi)= |D|_\mathfrak p^{-\frac{s+1}4}  L_\mathfrak p(s+1, \chi_\mathfrak p)
\frac{W_{t, \mathfrak p}(s, \phi)}{\gamma(W_{\mathfrak p})}
$$
with
$$
W_{t, \mathfrak p}(s, \phi)= \int_{F_\mathfrak p} \omega(wn(b))\phi_\mathfrak p (0) |a(wn(b))|_\mathfrak p^s\psi_\mathfrak p(-t b) \, db,
$$
and similarly
$$
W_{t, \sigma_j}^*(\tau_j, s, \Phi_{\sigma_j})= v_j^{-1/2} \pi^{-\frac{s+2}2} \Gamma(\frac{s+2}2) \int_{\R} \Phi_{\sigma_j}(w n(b)g_{\tau_j}, s) \psi(-bt) db.
$$
Then we have
\begin{equation}
E_t^*(\vec\tau, s, \phi)= \prod_{j=1}^2W_{t, \sigma_j}^*(\tau_j, s, \Phi_{\sigma_j}) \prod_{\mathfrak p\nmid \infty} W_{t, \mathfrak p}^*(s, \phi).
\end{equation}
Recall also (\cite[Proposition 2.6]{KRY1} that
$$
W_{\frac{t}{\sqrt D} , \sigma_j}^*(\tau_j, 0, \Phi_{\sigma_j}) =2i.
$$
Let  $\Diff(W, t/\sqrt D)$ be the set of prime ideals $\mathfrak p$  of $F$ such that $W_{\mathfrak p}$ does not represent $t/\sqrt D$. Then
$\mathfrak p  \in \Diff(W, t/\sqrt D)$ if and only if $ t \ne  2 u \bar u$ for any  $u \in {E_{\mathfrak p}}^\times$, i.e., $\mathfrak p$ is inert in $E/F$ and $\ord_{\mathfrak p}(2t)$ is odd (recall that $E/F$ is unramified), i.e.
$$
\Diff(W, t/\sqrt D) =\{ \mathfrak p \hbox{ inert in } E/F:\,  \ord_{\mathfrak p}(2t) \equiv 1 \pmod 2\}.
$$
It is a finite set of odd order. It is well-known that whenever $\mathfrak p \in  \Diff(W, t/\sqrt D)$, one has $W_{t, \mathfrak p}(0, \phi)=0$ for all $\phi \in S(W_\mathfrak p)$  (see for example \cite[Lemma 2.7]{YYCM}). So $a(t/\sqrt D, \phi) =0$  unless $\Diff(W, t/\sqrt D)=\{\mathfrak p\}$ has exactly one prime.  In this case,  we have
\begin{equation} \label{Firsteq}
a(\frac{t}{\sqrt D}, \phi_{\mu_5})
 =-4 W_{\frac{t}{\sqrt D}, \mathfrak p}^{*, \prime} (0,\phi_{\mu_5}) \prod_{\mathfrak q \nmid \mathfrak p \infty}
  W_{\frac{t}{\sqrt D}, \mathfrak q}^*(0, \phi_{\mu_5}).
\end{equation}

Let $\psi'(x) =\psi_F(\frac{x}{\sqrt D})$, which is unramified at every prime ideal of $F$, and $W'=W$ with $F$-quadratic form $Q_F'(x) =2 x\bar x$. Then the Weil representations associated to $(W, Q_F, \psi_F)$ and $(W', Q', \psi')$ are the same, and
$$
W_{\frac{t}{\sqrt D}, \mathfrak p}^{\psi_F}(s, \phi) = |D|_{\mathfrak p}^{\frac{1}4} W_{t, \mathfrak p}^{\psi'}(s, \phi)
$$
by (\ref{eq:ChangPsi}). So we have  in our case
\begin{equation} \label{eq4.7}
W_{\frac{t}{\sqrt D}, \mathfrak p}^*(s, \phi)=W_{t, \mathfrak p}^{*, \psi'}(s, \phi)
:=
|D|_{\mathfrak p}^{-\frac{s}4} L(s+1, \chi_{\mathfrak p})
  \frac{W_{t, \mathfrak p}^{\psi'}(s, \phi)}{\gamma(W_{\mathfrak p}')}.
  \end{equation}
Notice that $(\phi_{\mu_5}, \psi')$ is unramified at every prime  $\mathfrak p\nmid 2$, and we have by Proposition \ref{prop7.6}
$$
W_{\frac{t}{\sqrt D}, \mathfrak p}^*(0, \phi_{\mu_5})=\rho_{\mathfrak p}(t)
= \begin{cases}
  1+ \ord_{\mathfrak p}(t) &\ff E/F \hbox{ split  at } \mathfrak p,
  \\
  \frac{1+ (-1)^{\ord_\mathfrak p(t)}}2 &\ff E/F \hbox{ inert at } \mathfrak p.
  \end{cases}
  $$
In  particular,  $W_{\frac{t}{\sqrt D}, \mathfrak p}^*(0, \phi_{\mu_5})=0$  if and only if $\mathfrak  p \in \Diff(W', t)=\Diff(W, t/\sqrt D)$. In such a case, one has by the same proposition
$$
W_{\frac{t}{\sqrt D}, \mathfrak p}^{*, \prime}(0, \phi_{\mu_5}) =\frac{1 + \ord_{\mathfrak p}(t)}2.
$$
So we have proved the following formula: when  $\Diff(W, t/\sqrt D) = \{ \mathfrak p \}$, we have
\begin{equation}  \label{eq4.6}
a(\frac{t}{\sqrt D}, \phi_{\mu_5})
=-4 \begin{cases}
    \frac{1 +\ord_{\mathfrak p} (t\OO_F)}2\log\rm{N}(\frakp) \prod_{\mathfrak q\nmid 2 \mathfrak p} \rho_{\mathfrak q}(t)
      \prod_{\mathfrak q |2} W_{t, \mathfrak q}^{*,\psi'} (0, \phi_{\mu_5})  &\ff \mathfrak  p \nmid 2,
      \\
      W_{t, \mathfrak p}^{*,\psi', \prime}(0, \phi_{\mu_5}) \prod_{\mathfrak q\nmid 2 } \rho_{\mathfrak q}(t)  \prod_{\mathfrak q |2, \mathfrak q \ne \mathfrak p}
      W_{t, \mathfrak q}^{*, \psi'}(0, \phi_{\mu_5}) &\ff \mathfrak p|2.
      \end{cases}
\end{equation}

Now we focus on the local calculation at $\mathfrak p |2$.

\begin{lemma} Let $t =\frac{m+\sqrt D}4 \in \frac{1}2 \mathcal O_F-\mathcal O_F$ be as above, and let $\mathfrak p$ be  a prime of $F$ above $2$. Then $W_{t, \mathfrak p}^{*, \psi'}(s, \phi_{\mu_5}) =0$ unless $\alpha(\mu_5, t) =2 \mu_5 \bar{\mu}_5 -t \in  \OO_{\mathfrak p}$. Assuming $\alpha(\mu_5, t) \in \OO_{\mathfrak p}$ , we have the following.
\begin{enumerate}

\item When $F/Q$ is inert at $2$ with prime  $\mathfrak p$  over $2$.  Then $\mathfrak p \notin \Diff(W', t)$, $\ord_{\mathfrak  p} (2t) =0$,  and
$$
 W_{t, \mathfrak p}^{*, \psi'}(0, \phi_{\mu_5})  =\frac{1}4 L_{\mathfrak p}(1, \chi_{E/F}) =\frac{1}4 L_{\mathfrak p}(1, \chi_{E/F})  \rho_{\mathfrak p}(2t).
 $$

\item When $F/Q$ is ramified  at $2$ with $2\OO_F =\mathfrak  p^2$.  Then $\mathfrak p \notin \Diff(W', t)$, and
$$
 W_{t, \mathfrak p}^{*, \psi'}(0, \phi_{\mu_5})  =\frac{1}4 L_{\mathfrak p}(1, \chi_{E/F}) .
 $$
Moreover, assuming  $d_1 \equiv 0 \pmod 4$, one has
$$
\rho_{\mathfrak p}(2t) = \begin{cases}
  2 &\ff  d_2 \equiv 1 \pmod 8,
  \\
   1 &\ff d_2 \equiv 5 \pmod 8.
   \end{cases}
$$

\item
 When  $F/\Q$ is split at $2$,  let $\delta$ be the square root of $D$ in $\Z_2$ with $\delta \equiv 1 \pmod 4$, and let $2\OO_F=\mathfrak p_t \tilde{\mathfrak p}_t $ with $2t \in \mathfrak p_t$ and $2t \notin \tilde{\mathfrak p}_t $. Then
 \begin{enumerate}
 \item one has $\tilde{\mathfrak p}_t \notin \Diff(W', t)$, $\ord_{\tilde{\mathfrak  p}_t}(2t)=0$,  and
 $$
  W_{t, \tilde{\mathfrak p}_t}^{*, \psi'} (0, \phi_{\mu_5})  =\frac{1}2 L_{\tilde{\mathfrak p}_t}(1, \chi_{E/F}) =\frac{1}2 L_{\tilde{\mathfrak p}_t}(1, \chi_{E/F}) \rho_{\tilde{\mathfrak p}_t}(2t).
 $$

 \item
 $$
 W_{t,{\mathfrak p}_t}^{*, \psi'} (0, \phi_{\mu_5})  = \rho_{\mathfrak p_t}(t/2)=\rho_{\mathfrak p_t}(2t\mathfrak  p_t^{-2}).
 $$
In  particular, the   value is zero if and only if $\ord_{\mathfrak p_t} (2t)$ is odd and $E/F$ is inert at $\mathfrak p_t$ (i.e. $\mathfrak p_t \in \Diff(W', t)$). In this case, one has
$$
  W_{t,{\mathfrak p}_t}^{*, \psi', \prime} (0, \phi_{\mu_5})=\frac{1}2 \left(1+\ord_{\mathfrak p_t} (t/2)+L_{\mathfrak p_t}(1, \chi_{E/F})\right)\log \norm(\mathfrak p_t).
 $$
\end{enumerate}
\end{enumerate}
\end{lemma}

\begin{proof} We use the notation in Section \ref{sect5} and shorten $\mu =\mu_5 =\frac{-d_1 + d_2+\sqrt{d_1} +\sqrt{d_2}}4$. In particular,
$$
\alpha(\mu, t)  =-\frac{m+D}4 +\frac{d_1^2 -d_1 +d_2^2 -d_2}{8} \in  \OO_{\mathfrak p},
$$
and so $a(\mu, t) = \ord_{\mathfrak p}(\alpha(\mu, t) ) \ge 0$.

{\bf Case 1}: We first assume that $F/\Q$ is non-split at $2$ with  $\mathfrak p$ being the unique prime ideal of $F$ over $2$.  Since $t \notin \OO_{\mathfrak p}$ and $ \alpha (\mu, t) \in \mathcal O_{\mathfrak p}$, we have
$2 \mu \bar{\mu} \notin  \OO_{\mathfrak  p}$ and thus
$\ord_{E_{\mathfrak p}} (2 \mu) =0$,
where  for $z \in E_{\mathfrak p}$
$$
\ord_{E_{\mathfrak p}} (z) =\begin{cases}
   \ord_{\mathfrak p}(z) &\ff \mathfrak  p \hbox{ is inert in } E,
   \\
    \min (\ord_{\mathfrak B}(z),  \ord_{\bar{\mathfrak B}}(z) ) &\ff \mathfrak  p \OO_E = \mathfrak B \bar{\mathfrak  B} \hbox{ is split in } E.
\end{cases}
$$
So $a(\mu, t) \ge \ord_{ E_{\mathfrak p}}(2 \mu)$ and Proposition \ref{prop7.3} implies
$$
\frac{W_{t, \mathfrak p}^{\psi'}(0,  \phi_\mu)}{\gamma(W', \psi') \vol(\OO_{E_\mathfrak p}, d_{W'}x)} =1.
$$
Here $d_{W'}x$ is the self-dual Haar measure on $E_{\mathfrak p}=W'$ with respect to $(x, y) \mapsto \psi'(2 x \bar y)$, and so the volume
$$
\vol(\OO_{E_\mathfrak p}, d_{W'}x) = 1/4.
$$
Therefore,
$$
W_{t, \mathfrak p}^{*,\psi'}(0,  \psi_\mu) = \frac{1}4 L_\mathfrak  p(1,\chi_{E/F}).
$$
as claimed. Now we check
$$
\rho_{\mathfrak p}(2t)
  =\begin{cases}
    2 &\ff  4|d_i, d_j \equiv 1 \pmod 8,  i+j =3.
    \\
    1 &\hbox{ otherwise}.
    \end{cases}
$$
Assuming $4|d_1$ (the case $4|d_2$ is the same) first. We have  (notice that $D \equiv d_1 \pmod 8$)
$$
\norm_{E/F}(2\mu) = \frac{-D +\sqrt D}2 + \frac{d_1^2-d_1}4 + \frac{d_2^2 -d_2}4,
$$
and so
$$
\ord_{\mathfrak p} (\norm_{E/F}(2\mu)) =\begin{cases}
  1 &\ff d_2 \equiv 1 \pmod 8,
  \\
  0  &\ff d_2 \equiv 5 \pmod 8.
  \end{cases}
$$
Since
$$
2t - \norm_{E/F}(2\mu) = -2 \alpha(\mu, t) \in 2 \OO_F =\mathfrak p^2,
$$
we have that  $\ord_{\mathfrak p}(2t)$ is $1$ or $0$ depending whether $d_2 \equiv 1 \pmod 8$ or $5  \pmod 8$. This proves the ramified case.

Next assume that $F/\Q$ is inert at $2$. In this case, we have  (notice $\mathfrak p =2 \OO_F$)
$$
\ord_{\mathfrak  p} (\norm_{E/F}(2\mu)) = \ord_{\mathfrak p} (1+ \frac{-D+\sqrt D}2) =0,
$$
and thus $\ord_{\mathfrak p}(2t) =0$, $\rho_{\mathfrak p}(2t) =1$.
So we have proved (1) and (2) of the lemma.

{\bf Case 2}: We now assume that $F/\Q$ is split at $2$,  i.e. $D \equiv 1 \pmod 8$. This implies $d_1 \equiv d_2 \equiv 1, 5 \pmod 8$. Let $\delta \equiv 1 \pmod 4$ be one square root of $D$ in  $\Z_2$. Let $\mathfrak p_t$ be the prime of $F$ above $2$ such that $(\sqrt D)_{\mathfrak p_t} =\delta$, so $(\sqrt D)_{\tilde{\mathfrak p}_t}=-\delta$. It is easy to check (  recall $\sqrt{d_1} \sqrt{d_2} =-\sqrt D$)
\begin{align*}
\mu =(\mu_{\mathfrak p_t}, \mu_{\tilde{\mathfrak p}_t})  &\in   \OO_{E_{\mathfrak p_t}} \times (\frac{1}2\OO_{E_{\tilde{\mathfrak p}_t}} -\OO_{E_{\tilde{\mathfrak p}_t}}),
\\
 2 \mu \bar{\mu} &\in \OO_{\mathfrak p_t} \times (\frac{1}2\OO_{{\tilde{\mathfrak p}_t}} -\OO_{{\tilde{\mathfrak p}_t}}),
 \\
 t &\in \OO_{\mathfrak p_t} \times (\frac{1}2\OO_{{\tilde{\mathfrak p}_t}} -\OO_{{\tilde{\mathfrak p}_t}}),
\end{align*}
as  $\alpha(\mu, t)   \in \OO_F$.
 So Proposition \ref{prop7.3} asserts ($\beta =2$ and $\vol(\OO_{E_{{\mathfrak p}_t}},d_{\beta}x)=\vol(\OO_{E_{{\tilde{\mathfrak p}_t}}},d_{\beta}x)=1/2$ in this case)
$$
2\frac{W_{t, \tilde{\mathfrak p}_t}^{\psi'} (0, \phi_\mu)}{\gamma(W')} =1,
$$
i.e.,
$$
W_{t, \tilde{\mathfrak p}_t}^{*, \psi'}(0, \phi_\mu) = \frac{1}2 L_{\tilde{\mathfrak p}_t}(1, \chi_{E/F})=\frac{1}2 L_{\tilde{\mathfrak p}_t}(1, \chi_{E/F}) \rho_{\tilde{\mathfrak p}_t}(2t)),
$$
as $\ord_{\tilde{\mathfrak p}_t}(2t)=0$.

On the other hand,  Proposition \ref{prop7.6}(2) gives when $t \in  \mathfrak p_t $
$$
2L_{\mathfrak p_t}(1, \chi_{E/F}) \frac{W_{t, \mathfrak p_t}^{\psi'}(s, \phi_\mu)} {\gamma(W')}
=L_{\mathfrak p_t}(1, \chi_{E/F})(1-2^{-s}) + 2\sum_{0 \le n \le \ord_{\mathfrak p_t}(t/2)} (\chi_{\mathfrak p_t}(\pi_t) 2^{-s})^n.
$$
Here $\pi_t$ is a uniformizer of $F_{\mathfrak p_t}$. In  particular,
$$
W_{t, \mathfrak p_t}^{*, \psi'}(0, \phi_\mu) =\rho_{\mathfrak p_t} (t/2)
 =\rho_{\mathfrak p_t}(2t \mathfrak p_t^{-2}).
$$
Moreover, the value is zero if and only if $\mathfrak p_t$ is inert in $E/F$ (i.e., $d_i \equiv 5 \pmod 8$), and $\ord_{\mathfrak p_t} (2t)$ is odd, i.e., $\mathfrak p_t \in \Diff(W', t)$.  In such a case,
$$
W_{t, \mathfrak p_t}^{*, \psi', \prime}(0, \phi_\mu)=\frac{1}2 \left(1+\ord_{\mathfrak p_t} (t/2)+L_{\mathfrak p_t}(1, \chi_{E/F})\right)\log \norm(\mathfrak p_t).
$$
\end{proof}

In summary, we have the following proposition.

 \begin{proposition} \label{prop4.9} Let $t = \frac{m+\sqrt D}4 \in \frac{1}2 \OO_F$ with  $|m|< \sqrt D$.  Then the
 following is true.

\begin{enumerate}

\item  When  $|\Diff(W, \frac{t}{\sqrt D} ) | >1$,  we have $a(\frac{t}{\sqrt D},  \phi_{\mu_5}) =0$.

\item  When  $\Diff(W, \frac{t}{\sqrt D} ) =\{\mathfrak p \}$ with
$\mathfrak p \nmid 2$, let $\mathfrak p_2$ be the unique prime of $F$ over $2$. Then  we have
$$
a(\frac{t}{\sqrt{D}}, \phi_{\mu_5}) =-\epsilon(D) L_{\mathfrak p_2}(1, \chi_{E/F}) \frac{1+ \ord_{\mathfrak p}(t
\OO_F)}2
  \rho_{E/F}(2t\mathfrak p )\log\norm(\mathfrak p).
$$
Here
$$
\epsilon(D) = \begin{cases}
  2 &\ff D \equiv 1 \pmod 8,
  \\
  1 &\ff D \equiv  5 \pmod 8,
  \\
  \frac{1}2 &\ff D \equiv 0 \pmod 4.
  \end{cases}
$$

\item  When  $\Diff(W, \frac{t}{\sqrt D} ) =\{\mathfrak p \}$ with
$\mathfrak p | 2$, we have $d_1 \equiv d_2 \equiv 5 \pmod 8$, and  $2 \OO_F =\mathfrak  p \tilde{\mathfrak p}$
where $\ord_{\mathfrak p }(2t) >0$ is  odd and $\ord_{\tilde{\mathfrak p}}(2t) =0$.  In this case, we have
$$
a(\frac{t}{\sqrt{D}}, \phi_{\mu_5}) =-\epsilon(D) L_{\tilde{\mathfrak p}}(1, \chi_{E/F})
 \left[\frac{1+ \ord_{\mathfrak p}(2t \mathfrak p^{-1} )}2 +  \frac{1}3\right]
  \rho_{E/F}(2t\mathfrak p )\log\norm(\mathfrak p),
$$
with  $\epsilon(D) =2$ being as in  (2).
\end{enumerate}
 \end{proposition}

Here is a variant of the proposition, which will be used to prove the main theorem.

\begin{proposition}  \label{prop4.7} Let $t = \frac{m+\sqrt D}4 \in \frac{1}2 \OO_F$ with  $|m|< \sqrt D$.  Then the following is true.

\begin{enumerate}
\item When  $F/\Q$ is split at $2$ (i.e., $D \equiv 1 \pmod 8$), write $2\OO_F= \mathfrak p_t\tilde{\mathfrak p}_t$ with $t \in \OO_{F, \mathfrak p_t}$ (so $t \notin \OO_{F, \tilde{\mathfrak p}_t}$), one has
\begin{align*}
a(\frac{t}{\sqrt{D}}, \phi_{\mu_5}) &=-\epsilon(D)L_{\tilde{\mathfrak p}_t}(1, \chi_{E/F})
 \\
    & \quad \cdot \sum_{\mathfrak p  \hbox{ inert in } E/F}\left(\frac{1+ \ord_{\mathfrak p}(2 t\mathfrak p_t^{-2})}2+\frac{\ord_{\frakp_t}(\frakp)}{3}\right)\rho_{E/F}(2t\mathfrak p \mathfrak p_t^{-2})\log\norm(\mathfrak p).
\end{align*}

\item  When  $F/\Q$ is non-split at $2$, let $\mathfrak p_t =\OO_{F,2}$ and $\tilde{\mathfrak p}_t $ be the unique prime of $F$ over $2$. Then
$$
a(\frac{t}{\sqrt{D}}, \phi_{\mu_5}) =-\epsilon(D)L_{\tilde{\mathfrak p}_t}(1, \chi_{E/F})\sum_{\mathfrak p  \hbox{ inert in } E/F} \frac{1+ \ord_{\mathfrak p}(2t  \mathfrak p_t^{-2})}2
  \rho_{E/F}(2t\mathfrak p \mathfrak p_t^{-2})\log\norm(\mathfrak p).
$$
\end{enumerate}

\end{proposition}

\subsection{Computing $a_0(\phi_\mu)$}

Let $\psi'(x) =\psi_F(\frac{x}{\sqrt D})$ and  $W'=E$ with $Q_F'(x)=2 x \bar x$ as before. Proposition 5.7 (3) and Remark \ref{rem:Constant} imply
\begin{equation} \label{eq4.8}
a_0(\phi_\mu)=-\Lambda(0,\chi)\frac{d}{ds} \prod_{\mathfrak p|2} \tilde{W}_{0, \mathfrak p}^{\psi'} (0, \phi_\mu)|_{s=0},
\end{equation}
where
\begin{equation}  \label{eq4.9}
\tilde{W}_{0, \mathfrak p}^{\psi'}(s,\phi_{\mu}) =\frac{L_\mathfrak p(s+1, \chi)}{L_\mathfrak p(s, \chi)} |2|_\mathfrak p \frac{W_{0, \mathfrak p}^{\psi'}(s, \phi_\mu)}{\gamma(W_{\mathfrak p}') \vol(\OO_{E_\mathfrak p}, d_{Q_F'}x, \psi')}.
\end{equation}
By Corollary \ref{cor:ConstantTerm}, one has $a_0(\phi_{\mu_6})=0$. For other $\mu_j$s, it is complicated and we need to do it on  a case-by-case basis. The main tools are Corollary \ref{cor:ConstantTerm} and Proposition \ref{prop7.6}.
Let
\begin{equation}
b_i=\ord_2(\frac{d_i^2-d_i}4)=\ord_2(e_{i+1}\bar{e}_{i+1})=
\begin{cases}
 0 &\ff d_i \equiv 4, 5\pmod 8,
 \\
 1 &\ff d_i \equiv 0 \pmod 8,
 \\
  \ord_2(d_i-1)  -2 \ge 1 &\ff d_i \equiv 1 \pmod 8,
\end{cases}
\end{equation}
and
$$
 c_1=\ord_2( (1+e_{2})(1+\bar{e}_{2}))=\ord_2 (1-d_1 + \frac{d_1^2 -d_1}4)
  =\begin{cases}
    0 &\ff d_1 \equiv 0, 5 \pmod 8,
    \\
    1 &\ff d_1 \equiv 4 \pmod 8,
    \\
    \ge 2 &\ff d_1 \equiv 1 \pmod {16},
    \\
     1 &\ff d_1 \equiv 9 \pmod {16}.
     \end{cases}
$$
$$
 c_2=\ord_2( (1+e_{3})(1+\bar{e}_{3}))=\ord_2 (1+d_2 + \frac{d_2^2 -d_2}4)
  =\begin{cases}
    0 &\ff d_2 \equiv 0, 5 \pmod 8,
    \\
    1 &\ff d_2 \equiv 4 \pmod 8,
    \\
    1 &\ff d_2 \equiv 1 \pmod {16},
    \\
     \ge 2 &\ff d_2 \equiv 9 \pmod {16}.
     \end{cases}
$$
 We denote $\mathfrak p$ (or $\mathfrak p_i$) for primes of $F$ above $2$, and $\mathfrak B$ and $\bar{\mathfrak B}$ (or $\mathfrak B_i$ and  $\bar{\mathfrak B}_i$) for primes of $E$ above $\mathfrak p$ (or $\mathfrak p_i$). First we have the following table, where $a_\mathfrak p(\mu, 0)=\ord_\mathfrak p(2 \mu \bar\mu)$.

\begin{table}[h]
\caption{\label{table2} When does $a_{\mathfrak p}(\mu,0 )=\ord_{\mathfrak p}(2 \mu \bar\mu)  \ge 0$?}
\begin{tabular}{|c|c|c|c| }
\hline \rule[-3mm]{0mm}{8mm}
$d_1 \pmod 8$ & $d_2 \pmod 8$ & $a_{\mathfrak p}(\mu, 0) <0$              &$a_{\mathfrak p}(\mu, 0) \ge 0$  \\
\hline  $1$   &$0$         & $ \mu_7$                                     &$\mu_1, \mu_2, \mu_3,\mu_4,\mu_9, \mu_{10}, \mu_{12}$ \\
\hline  $1$   &$1$         &                                              &$\mu_1, \mu_2, \mu_3,\mu_4,\mu_7,\mu_9, \mu_{10}, \mu_{12}$ \\
\hline  $1$   &$4$         & $ \mu_4$                                     &$\mu_1,\mu_2, \mu_3,\mu_7,\mu_9, \mu_{10}, \mu_{12}$ \\
\hline  $1$   &$5$         & $ \mu_4, \mu_7$                              &$\mu_1, \mu_2, \mu_3,\mu_9, \mu_{10}, \mu_{12}$ \\
\hline  $5$   &$0$         &$\mu_1, \mu_3,\mu_7,\mu_9,\mu_{10}$     &$\mu_2, \mu_4, \mu_{12}$\\
\hline  $5$   &$4$         & $\mu_2, \mu_3,\mu_4,\mu_{10}, \mu_{12}$      &$\mu_1,\mu_7,\mu_9$ \\
\hline  $5$   &$5$         &$\mu_1,\mu_2, \mu_3,\mu_4, \mu_7,\mu_9, \mu_{10},\mu_{12}$       &     \\
\hline
\end{tabular}
\end{table}

\begin{lemma}\label{d5} Recall that $(d_1, d_2) =1$.  Then the following is true.
\begin{enumerate}
\item When $d_1 \equiv d_2 \equiv 5 \pmod 8$, one has $a_0(\phi_\mu)=0$ for $\mu=\mu_1,\mu_2,\mu_3, \mu_4, \mu_7,\mu_9, \mu_{10},\mu_{12}$.

\item  When  $d_1 \equiv 5 \pmod 8$, and $d_2 \equiv 4 \pmod 8$, we have
$$
a_0(\phi_\mu) =\begin{cases}
   0 &\ff \mu_2,  \mu_3, \mu_{4}, \mu_{10}, \mu_{12},
   \\
    -\frac{\Lambda(0,\chi)}3 \log 2 &\ff \mu=\mu_1, \mu_7,\mu_9.
    \end{cases}
$$

\item When  $d_1 \equiv 5 \pmod 8$ and $d_2 \equiv 0 \pmod 8$, we have
$$
a_0(\phi_\mu) =\begin{cases}
   0 &\ff \mu_1, \mu_3,  \mu_7, \mu_{9}, \mu_{10},
   \\
    -\frac{\Lambda(0,\chi)}3 \log 2 &\ff \mu=\mu_2,\mu_4,\mu_{12}.
    \end{cases}
$$

\end{enumerate}
\end{lemma}
\begin{proof} Claim (1) follows from  Table \ref{table2} and Corollary \ref{cor:ConstantTerm}.  For (2) and (3), let $\mathfrak p$ be the unique ramified prime of $F$ above $2$ and let $\mathfrak B =\mathfrak p \OO_E$ be the unique inert prime of $E$ above $\mathfrak p$.

When $d_2 \equiv 4 \pmod 8$, one has $d_2 \equiv -4 \pmod {16}$ (as $d_2/4 \equiv -1 \pmod 4$) and so $c_2 = 1$.
As $(1+e_3) (1+\bar{e}_3) \in \Q_2$, we have
$$
\ord_{\mathfrak B}(\frac{1}2 (1+e_3)) =
 \frac{1}2 \ord_{\mathfrak p} (\frac{1}4 (1+e_3) (1+\bar{e}_3) ) =
 \ord_2 (\frac{1}4 (1+e_3) (1+\bar{e}_3) ) =
 (c_2 -2)=-1.
$$
So   $\mu  \notin \OO_{\mathfrak B}$ and $a_\mathfrak p(\mu, 0) =0$ for $\mu=\mu_1, \mu_7,\mu_9$.  Now (2) follows from Corollary \ref{cor:ConstantTerm}.

When $d_2 \equiv 0 \pmod 8$, one has $a_\mathfrak p(\mu, 0) =0$ for  $\mu=\mu_2, \mu_{4},\mu_{12}$.  So (3) follows from Corollary \ref{cor:ConstantTerm}.

\end{proof}

When  $d_1 \equiv 1 \pmod 8$, and $d_2 \not\equiv 1 \pmod 8$, let $\mathfrak p$ be the unique prime of $F$ over $2$, then $\mathfrak p=\mathfrak P \bar{\mathfrak P}$ is split in $E$. Let $\sqrt{d_1}$ be the unique square root of $d_1$ in $\Q_2$ with $\sqrt{d_1}\equiv 1 \pmod 4$, and identify
$$
E \subset  E_\mathfrak p =E_{\mathfrak P} \times E_{\bar{\mathfrak P}} \cong F_{\mathfrak p} \times F_{\mathfrak p}, \quad \mu \mapsto (\mu, \bar\mu).
$$
Let $o_{E_{\mathfrak p}}(\mu) = \min (o_{\mathfrak p}(\mu), o_{\mathfrak p}(\bar\mu))$ as in  (\ref{eq5.4}). Direct calculation  gives the following table:

\begin{table}[h]
\caption{\label{table3} $o_{E_{\mathfrak p}}(2\mu)$}
\begin{tabular}{|c|c|c|c| }
\hline \rule[-3mm]{0mm}{8mm}
                                         &$d_2 \equiv 5 \pmod 8$          &$d_2 \equiv 0 \pmod 8$ &$d_2 \equiv 4 \pmod 8$ \\
\hline  $\mu_1=\frac{1}2 e_2(1+e_3)$        &$0$                              &$0$                   &$1$\\
\hline  $\mu_2=\frac{1}2 e_2 e_3$        &$0$                              &$1$                   &$0$\\
\hline  $\mu_3=\frac{1}2 e_2 $        &$0$                              &$0$                   &$0$\\
\hline  $\mu_4=\frac{1}2 e_3$        &$0$                              &$1$                   &$0$\\
\hline  $\mu_7 =\frac{1}2(1+e_3)$        &           $0$                    &  $0$             &$1$\\
\hline  $\mu_9=\frac{1}2 (1+e_2) (1+e_3)$        &$0$                              &$0$                   &$1$\\
\hline $\mu_{10}=\frac{1}2 (1+e_2)$                           & $0$                          &$ 0$                    &$0$\\
\hline $\mu_{12}=\frac{1}2 e_3(1+e_2)$                           & $0$                          &$1$                    &$0$\\
\hline
\end{tabular}
\end{table}

Notice also that $a_{\mathfrak p}(\mu, 0) = \ord_{\mathfrak p}(2 \mu \bar\mu) <0$ for $\mu=\mu_6$ always and for $\mu=\mu_7$ with $d_2 \equiv 0,5 \pmod 8$ and for $\mu=\mu_4$ with $d_2 \equiv 4,5 \pmod 8$.
\begin{lemma} Assume  that $d_1 \equiv 1 \pmod 8$, and $d_2 \not\equiv 1 \pmod 8$. Then  $a_0(\phi_{\mu_6})=0$, and
\begin{enumerate}
\item  When  $d_2 \equiv 5 \pmod 8$, one has
$$
a_0(\phi_{\mu_j}) = \begin{cases}
   0 &\ff   j=4,7,
   \\
   -\frac{2\Lambda(0,\chi)}3 \log 2 &\ff j=1,2, 3, 9,10,12.
\end{cases}
$$

\item  When $d_2 \equiv 0 \pmod 8$, one has
$$
a_0(\phi_{\mu_j}) = \begin{cases}
   0 &\ff   j=4, 7,
   \\
   -\frac{\Lambda(0,\chi)}2 \log 2 &\ff j=1, 3, 9,10
   \\
  - \Lambda(0,\chi) \log 2 &\ff j=2,12.
\end{cases}
$$

\item When $d_2\equiv 4 \pmod 8$, one has
$$
a_0(\phi_{\mu_j}) = \begin{cases}

   0&\ff j=4,7\\
  - \frac{\Lambda(0,\chi)}2 \log 2 &\ff j=2, 3, 10, 12,
   \\
  - \Lambda(0,\chi) \log 2 &\ff j=1,9.
\end{cases}
$$

\end{enumerate}

\end{lemma}

\begin{proof}
Let $\mu=\mu_1, \mu_2,\mu_3,\mu_{9},\mu_{10},\mu_{12}$, then $a_\frakp(\mu,0)\geq 0$. We see that if $o_{E_{\mathfrak p}}(2\mu)=0$ then only one of $\mu$ and $\bar{\mu}$ belong to $\mathcal{O}_{F_{\frakp}}$. We can also check directly that if $o_{E_{\mathfrak p}}(2\mu)=1$ then
$a_\frakp(\mu,0)\geq 1$, and again only one of $\mu$ and $\bar{\mu}$ belong to $\mathcal{O}_{F_{\frakp}}$. Then using table \ref{table3}, equations (\ref{eq4.8}), (\ref{eq4.9}) and  Corollary \ref{cor:ConstantTerm} we can get  the formulae in the lemma by noticing that $2$ is inert in $F$ if $d\equiv 5\mod 8$ and ramified otherwise. The cases  $\mu=\mu_4, \mu_7$  can also be dealt with similarly by noting that both $\mu$ and $\bar{\mu}$ belong to $\mathcal{O}_{F_{\frakp}}$ in both cases.

\end{proof}

Finally, we deal with the case $d_1 \equiv d_2 \equiv 1 \pmod 8$.  In this case, we have
 the following lemma.

\begin{lemma}\label{d11} Assume $d_i \equiv 1 \pmod 8$ with $i=1, 2$. We have
$$
a_0(\phi_{\mu_j}) =\begin{cases}
0 &\ff  j=3,4,7,10,
\\
-\Lambda(0, \chi) \log 2 &\ff j=1, 2, 9, 12.
 \end{cases}
$$
\end{lemma}

\begin{proof}
 Write
$$
2\OO_F =\mathfrak p_1 \mathfrak p_2, \quad \mathfrak p_i \OO_E =\mathfrak P_i \bar{\mathfrak P}_i.
$$
Let  $\sqrt{d_i}  \equiv 1 \pmod 4 \in \Z_2$ be the prefixed square roots of $d_i$ and  let $ \sqrt D =-\sqrt{d_1} \sqrt{d_2} \in  \Z_2$. We identify $F_{\mathfrak p_i}$, $E_{\mathfrak P_i}$, and $E_{\bar{\mathfrak P}_i}$ with $\Q_2$ as follows:
\begin{align*}
F_{\mathfrak p_i} &\cong \Q_2, \quad \sqrt D \mapsto (-1)^{i-1}\sqrt D,
\\
E_{\mathfrak P_i}&\cong \Q_2,  \quad \sqrt D \mapsto (-1)^{i-1}\sqrt D, \sqrt{d_i} \mapsto \sqrt{d_i},
\\
 E_{\bar{\mathfrak P}_i}&\cong \Q_2,  \quad \sqrt D \mapsto (-1)^{i-1}\sqrt D, \sqrt{d_i} \mapsto -\sqrt{d_i}.
\end{align*}
We also fix the following embedding
$$
E \subset E\otimes_\Q \Q_2 = E_{{\mathfrak P}_1} \times E_{\bar{\mathfrak P}_1} \times E_{{\mathfrak P}_2} \times E_{\bar{\mathfrak P}_2} \cong \Q_2^4
$$
such that
$$
e_2 \mapsto (e_2, \bar{e}_2, \bar{e}_2, e_2), \quad e_3 \mapsto (e_3, \bar{e}_3, e_3, \bar{e}_3).
$$
With this identification and direct calculation, we have Table \ref{table4}.

\begin{table}[h]
\caption{ Valuations of $\mu$--- when $d_1, d_2 \equiv 1 \pmod 8$}\label{table4}
\begin{tabular}{|c|c|c|c|c|c|c| }
\hline \rule[-3mm]{0mm}{8mm}
    $\mu$  &$\ord_{\mathfrak P_1}\mu$     &$\ord_{\bar{\mathfrak P}_1}\mu$     &$\ord_{\mathfrak P_2}\mu$   &$\ord_{\bar{\mathfrak P}_2}\mu$
     & $\mu \in \mathcal O_{E_{\mathfrak p_1}} ?$  & $\mu \in \mathcal O_{E_{\mathfrak p_2}} ?$  \\
\hline $e_2$                            &$\ge 1$          &  $0$         &$0$         &$\ge 1$        &Yes  &Yes\\
\hline  $1+e_2$                         &$0$             &$\ge 1$      &$\ge 1$      &$0$            &Yes &Yes\\
\hline $e_3 $                           &$0$              &  $\ge 1$    &$0 $          & $\ge 1$           &Yes  &Yes\\
\hline  $1+e_3$                         &$\ge 1$          &$0$           &$\ge 1$     &$0$          &Yes  &Yes  \\
\hline $\mu_1= \frac{1}2 e_2(1+ e_3)$       &$\ge 1$         &$-1$      & $\geq 0$       &$\geq 0$        &No   &Yes  \\
\hline $\mu_2= \frac{1}2 e_2 e_3$       &$\ge 0$         &$\ge 0$      & $-1$       &$\ge 0$        &Yes   &No  \\
\hline $\mu_3= \frac{1}2 e_2$           &$\ge 0$         &$-1$      & $-1$       &$\ge 0$        &No   &No  \\
\hline $\mu_4= \frac{1}2 e_3$           &$-1$            &$\ge 0$      & $-1$       &$\ge 0$        &No   &No  \\
\hline $\mu_6= \frac{1}2$               &$-1$            &$-1$      & $-1$       &$-1$        &No   &No  \\
\hline $\mu_7= \frac{1}2 (1+e_3)$       &$\ge 0$         &$-1$      & $\ge 0$       &$-1$        &No  &No  \\
\hline $\mu_9= \frac{1}2 (1+e_2)(1+e_3)$   &$\ge 0$      &$\ge 0$      & $\ge 0$       &$-1$        &Yes   &No  \\
\hline $\mu_{10}= \frac{1}2 (1+e_2)$       &$-1$         &$\ge 0$      & $\ge 0$       &$-1$        &No   &No  \\
\hline $\mu_{12}= \frac{1}2 e_3(1+e_2)$       &$-1$         &$\ge 1$      & $\ge 0$       &$\geq 0$        &No   &Yes  \\
\hline
\end{tabular}
\end{table}

By (\ref{eq4.8}) and (\ref{eq4.9}), we have
$$
a_0(\phi_\mu) =-\Lambda(0, \chi) \frac{d}{ds} \left[\prod_{j=1}^2 \frac{W_{0, \mathfrak p_j}^{\psi'}(s, \phi_\mu)}{\gamma(W_{\mathfrak p_j}') \vol (\OO_{E_{\mathfrak p_j}}, d_2 x)}\right]_{s=0}
$$
with
$$
\tilde{W}_{0, \mathfrak p_j}(s, \phi_\mu) = \frac{1}2 \frac{1-2^{-s}}{1-2^{-1-s}} \frac{W_{0, \mathfrak p_j}^{\psi'}(s, \phi_\mu)}{\gamma(W_{\mathfrak p_j}') \vol (\OO_{E_{\mathfrak p_j}}, d_2 x)}.
$$
When  $\mu \notin \mathcal O_{E_{\mathfrak P_j}}$, $W_{0, \mathfrak p_j}^{\psi'}(s, \phi_\mu)$ is a polynomial of $2^{-s}$, and so
$$
\tilde{W}_{0, \mathfrak p_j}(0, \phi_\mu)=0
$$
and
$$
\tilde{W}_{0, \mathfrak p_j}'(0, \phi_\mu)=\frac{W_{0, \mathfrak p_j}^{\psi'}(0, \phi_\mu)}{\gamma(W_{\mathfrak p_j}') \vol (\OO_{E_{\mathfrak p_j}}, d_2 x)} \log 2.
$$
This implies that $a_0(\phi_\mu)=0$ for $\mu =\mu_3,\mu_4,\mu_7, \mu_{10}$.

When  $\mu \in \mathcal O_{E_{\mathfrak P_j}}$, $\phi_\mu =\phi_0$,  one has by Proposition \ref{prop7.6}(3)
$$
\tilde{W}_{0, \mathfrak p_j}(s, \phi_\mu) =1 + \frac{1}2(1-2^{-s})\cdot (\hbox{ polynomial of } 2^{-s}).
$$
Thus $\tilde{W}_{0, \mathfrak p_j}(0, \phi_\mu)=1$.

For $\mu=\mu_2$, the table shows that $\mu  \in \OO_{E_{\mathfrak p_1}}$ and $\mu \notin \OO_{E_{\mathfrak p_2}}$, so we have
$$
a_0(\phi_\mu) =-\Lambda(0, \chi) \tilde{W}_{0, \mathfrak p_1} (0, \chi_{\mu_2}) \tilde{W}_{0, \mathfrak p_2}'(0, \phi_{\mu_2})
             =-\Lambda(0, \chi) \log 2.
$$
The cases for $\mu_1, \mu_{9},\mu_{12}$ are similar and left to the reader.

\end{proof}

In summary, we have the following lemma which is needed to prove Theorem \ref{maintheorem} and Proposition \ref{prop:norm}.

\begin{proposition}  \label{prop:Constant} The following are true.
\begin{enumerate}
\item  One has
\begin{align*}
 &8a_0(\phi_{\mu_2}) -8a_0(\phi_{\mu_3})-8a_0(\phi_{\mu_4})-8a_0(\phi_{\mu_6})-8a_0(\phi_{\mu_7}) +8a_0(\phi_{\mu_{9}})-8a_0(\phi_{\mu_{10}})
 \\
 &= -4 \Lambda(0,\chi) \log 2 \begin{cases}
 4 &\ff (d_1, d_2) \equiv (1 \mod 8, 1 \mod 8),
 \\
 1 &\ff  (d_1, d_2) \equiv (1 \mod 8, 0 \mod 4) \hbox{ or }  (0 \mod 4, 1 \mod 8),
 \\
  0   &\hbox{otherwise}.
  \end{cases}
\end{align*}

\item  One has
\begin{align*}
&a_0(\phi_{\mu_2})-a_0(\phi_{\mu_3})+a_0(\phi_{\mu_4})-a_0(\phi_{\mu_6})-a_0(\phi_{\mu_{10}})+a_0(\phi_{\mu_{12}})
\\
&=\Lambda(0, \chi) \log 2   \begin{cases}
 0  &\ff d_1 \equiv 5 \pmod 8,  d_2 \equiv  4, 5 \pmod 8,
 \\
 -1  &\ff d_1 \equiv 5 \pmod 8,  d_2 \equiv 0 \pmod 8,
 \\
  -2  &\ff d_1 \equiv d_2 \equiv 1 \pmod 8.
  \end{cases}
  \end{align*}

  \item  One has
  \begin{align*}
  &a_0(\phi_{\mu_1})-a_0(\phi_{\mu_3})-a_0(\phi_{\mu_6})+a_0(\phi_{\mu_7})+a_0(\phi_{\mu_9})-a_0(\phi_{\mu_{10}})
  \\
  &= \Lambda(0, \chi) \log 2 \begin{cases}
   0  &\ff d_1 \equiv 5 \pmod 8,  d_2 \equiv 0, 5\pmod 8,
 \\
 -1  &\ff d_1 \equiv 5 \pmod 8,  d_2 \equiv 4 \pmod 8,
 \\
  -2  &\ff d_1 \equiv d_2 \equiv 1 \pmod 8.
  \end{cases}
  \end{align*}

\end{enumerate}

\end{proposition}

{\bf Proof of Theorem  \ref{maintheorem}}: Now Theorem  \ref{maintheorem} follows from  Propositions \ref{prop:norm2}, \ref{prop4.7}, and \ref{prop:Constant}(1), and the fact that
$$
\Lambda(0, \chi) =\Lambda(0, \chi_{\kay_1/\Q}) \Lambda(0, \chi_{\kay_2/\Q}) = \frac{4 h_1 h_2}{w_1 w_2}.
$$

{\bf Proof of Proposition \ref{prop:norm}}: Take $d_2 =d$, $\tau_2=\tau =\frac{d+\sqrt d}2$. Choose some $d_1 <-4$  with   $(d_1, d_2) =1$,  then we have by Corollary \ref{lambda} and Theorem \ref{theo:BigCM} (recall $\lambda_0=\lambda(\frac{d+\sqrt d}2)$) that
 $$
 \log|\norm(\lambda_0)|
=-\frac{r(d_1)r(d_2)}{4\delta(d_2)h(\kay_1,2)}8(a_0(\phi_{\mu_2})-a_0(\phi_{\mu_3})+a_0(\phi_{\mu_4})-a_0(\phi_{\mu_6})-a_0(\phi_{\mu_{10}})+a_0(\phi_{\mu_{12}})),
 $$
Now applying Proposition \ref{prop:Constant} (2), we obtain via a little calculation
$$
\log|\norm(\lambda_0)|=\begin{cases}
 0 &\ff  d \equiv 4, 5 \pmod 8,
 \\
 2 h &\ff d \equiv 0 \pmod 8,
 \\
 4 h &\ff d \equiv 1 \pmod 8
 \end{cases}
$$
as claimed.  The same argument with
$$
\log|\norm(1-\lambda_0))|
=-\frac{r(d_1)r(d_2)}{4\delta(d_2)h(\kay_1, 2)}8(a_0(\phi_{\mu_1})-a_0(\phi_{\mu_3})-a_0(\phi_{\mu_6})+a_0(\phi_{\mu_7})+a_0(\phi_{\mu_9})-a_0(\phi_{\mu_{10}}))
$$
gives the desired  formula for  $|\norm(1-\lambda_0))|$. This proves Proposition \ref{prop:norm}.

\subsection{Proof of Theorem \ref{theo:Difference}} In  this subsection, we will study the set $\mathcal N(d_1, d_2)$ and prove Theorem \ref{theo:Difference}. For an  imaginary quadratic field $\kay =\Q(\sqrt d)$, set $\tau^0 =\frac{d+\sqrt{d}}2$. We also write
$\tau_1^0 =\frac{d_1+\sqrt{d_1}}2$ and $\tau_2^0 =\frac{d_2+\sqrt{d_2}}2$. Denote
$$
\lambda_1= \lambda,~\lambda_2=1-\lambda,~\lambda_3=\frac{1}{\lambda},~\lambda_4=\frac{\lambda}{\lambda-1},
~\lambda_5=\frac{\lambda-1}{\lambda}, \hbox{ and } ~\lambda_6=\frac{1}{1-\lambda}.
$$
Then the values $\lambda(\tau)$, $\tau \in \CM(\kay,2)$, are determined up to Galois conjugation by $\lambda_i(\tau^0)$,  $1\le i \le 6$.

\begin{lemma}
\label{lem:norm}
Let  $a(d) =h(\kay, 2)$ or $2h(\kay, 2)$ depending on whether $4|d$ or not. Then
\begin{align*}
\mathcal N(d_1, d_2)
 &= \cup \left\{
 \left| \norm(\lambda(\tau_1^0) - \lambda_i(\tau_2^0))\right|:
   \, 1 \le i \le 6 \right\}
\\
 & \cup \left\{
 \left| \frac{\norm(\lambda(\tau_1^0) - \lambda_i(\tau_2^0))}{\norm(\lambda(\tau_1^0)^{a(d_2)})\norm(\lambda_i(\tau_2^0)^{a(d_1)})}\right|:
   \, 1 \le i \le 6 \right\}
   \\
   & \cup \left\{ \left| \frac{\norm(\lambda(\tau_1^0) - \lambda_i(\tau_2^0))}{\norm((1-\lambda(\tau_1^0))^{a(d_2)})\norm((1 -\lambda_i(\tau_2^0))^{a(d_1)})} \right|:\, 1 \le i \le 6\right\}.
\end{align*}
\end{lemma}
\begin{proof} For $\tau_l \in \CM(\kay_l, 2)$ with  $l=1, 2$, we may assume (up to Galois conjugation) $\lambda(\tau_1) =\lambda_i(\tau_1^0)$ and $\lambda (\tau) =\lambda_i(\tau_2')$ for some $i$.
 Simple computation gives
$$\lambda(\tau_1)-\lambda(\tau_2)=\lambda_i(\tau_1^0)-\lambda_i(\tau'_2)=\pm
\left\{\begin{array}{ll}
\lambda(\tau_1^0)-\lambda(\tau'_2)&\text{if }i=1,2,\\
\frac{\lambda(\tau_1^0)-\lambda(\tau'_2)}{\lambda(\tau_1^0)\lambda(\tau'_2)}&\text{if }i=3,5,\\
\frac{\lambda(\tau_1^0)-\lambda(\tau'_2)}{(1-\lambda(\tau_1^0))(1-\lambda(\tau'_2))}&\text{if }i=4,6.
\end{array}\right.$$
Now write $\lambda(\tau_2') =\lambda_i(\tau_2^0) $ for some $i$. Taking norm on both sides of the above identity  proves  the lemma.
\end{proof}

 Let
$$
\mathcal N_0(d_1,d_2)=\{  | \norm(\lambda(\tau_1^0) - \lambda_i(\tau_2^0))|:\, 1 \le i \le 6\}
=
\{  | \norm(\lambda(\tau_1^0) - \lambda(\tau_2))|:\, \tau_2 \in  \CM(\kay_2, 2)\}.
$$
Then $|\mathcal N_0(d_1,d_2)|=1$ when  $d_2 \equiv 5 \pmod 8$ and $ |\mathcal N_0(d_1,d_2)| \le 3$ otherwise by Proposition \ref{prop:MininalPoly}.

\begin{lemma}
\label{lem:numerator}
Assume $d_1\equiv 1 \pmod 8$ and $d_2 \not\equiv 5 \pmod 8$.  Then $|\mathcal N_0(d_1, d_2)|=2$. More precisely, write
$ \mathcal N_0(d_1, d_2)=\{ \alpha,  \beta, \gamma\}$ with
\begin{eqnarray*}
\alpha&=&|\norm(\lambda(\tau_1^0) - \lambda(\tau_2^0))|=|\norm(\lambda(\tau_1^0) - \lambda_2(\tau_2^0))|,\\
\beta&=&\left|\norm\left(\lambda(\tau_1^0) - \lambda_3(\tau_2^0)\right)\right|=\left|\norm\left(\lambda(\tau_1^0) - \lambda_5(\tau_2^0)\right)\right|,\\
\gamma&=&\left|\norm\left(\lambda(\tau_1^0) - \lambda_4(\tau_2^0)\right)\right|=\left|\norm\left(\lambda(\tau_1^0) - \lambda_6(\tau_2^0)\right)\right|.
\end{eqnarray*}
 Then  the following are true.

  \begin{table}[h]
  \caption{\label{table5} $\alpha,\beta,\gamma$}
  \begin{tabular}{|c|c|c|c|}
  \hline \rule[-3mm]{0mm}{8mm}
  $d_2\pmod 8$ & $4$ & $0$ &$1$\\
  \hline
   & $\alpha =\beta$ & $\alpha =\gamma$ &$\beta =\gamma$ \\
  \hline
  \end{tabular}
  \end{table}

  \begin{table}[h]
  \caption{\label{table6}2-parts in $\alpha,\beta,\gamma$}
  \begin{tabular}{|c|c|c|}
  \hline \rule[-3mm]{0mm}{8mm}
  $d_2 \pmod 8$ & $4 $ \hbox{ or  } $0$  & $1$\\
  \hline
  $\ord_2(\alpha)$ & $2h_1h_2$ & $8h_1h_2$\\
  \hline
  $\ord_2(\tilde{\alpha})$ & $-4h_1h_2$ & $-4h_1h_2$\\
  \hline
  \end{tabular}
  \end{table}
  Here $\tilde{\alpha}$ is the element of  $\{\beta,\gamma\}$ that is different from $\alpha$.

\end{lemma}

\begin{proof} Proposition \ref{prop:GaloisConjugate} asserts that $\lambda(\tau_1^0)$ and $1-\lambda(\tau_1^0) $ are Galois conjugates of each other. So
$$
|\norm(\lambda(\tau_1^0)-\lambda_i(\tau_2^0)) |= |\norm(\lambda(\tau_1^0) - (1-\lambda_i(\tau_2^0))|,
$$
which justify the two formulas for $\alpha$, $\beta$ and $\gamma$. Table \ref{table5} also follows directly from  Proposition \ref{prop:GaloisConjugate}.

Notice that $\lambda_i(\tau_2^0)$,  $1 \le i \le 6$, are the six roots of
$f(\lambda,j(\tau_2^0))$  in equation (\ref{eq:generic}). So we have
$$f(\lambda,j(\tau_2^0))=\prod_{i=1}^6(\lambda-\lambda_i(\tau_2^0)).$$
In particular,
$$f(\lambda(\tau_1^0),j(\tau_2^0))=\prod_{i=1}^6(\lambda(\tau_1^0)-\lambda_i(\tau_2^0)).$$

On the other hand, $\lambda(\tau_1^0)$ is also a root of $f(\lambda,j(\tau_1^0))$.  Therefore,
\begin{eqnarray}
\label{eq:relation}
256\prod_{i=1}^6(\lambda(\tau_1^0)-\lambda_i(\tau_2^0))&=&256(f(\lambda(\tau_1^0),j(\tau_2^0))-f(\lambda(\tau_1^0),j(\tau_1^0)))\\
&=&(j(\tau_1^0)-j(\tau_2^0))\lambda(\tau_1^0)^2(1-\lambda(\tau_1^0))^2.\nonumber
\end{eqnarray}

Let $H(\kay_i,2)$ be the ray class field of $\Q(\sqrt{d_i})$ of modulus 2, $i=1,2$. Then we have $[H(\kay_i,2):\Q]=2h_i$ by Lemma \ref{lem:ClassGroup} and we also know that $H(\kay_1, 2)$ and $H(\kay_2, 2)$ are disjoint by Lemma \ref{lem:ClassIsomorphism}.
Taking norms in $H(\kay_1,2)H(\kay_2,2)/\Q$ on both sides of equation (\ref{eq:relation}), we have
$$256^{4h_1h_2}(\alpha\beta\gamma)^2=(\norm(j(\tau_1^0)-j(\tau_2^0)))^{4}\norm(\lambda(\tau_1^0))^{4h_2}\norm((1-\lambda(\tau_1^0)))^{4h_2}.$$
Applying Proposition \ref{prop:norm}, we obtain
\begin{equation} \label{eq:NormRelation}
\norm(j(\tau_1^0)-j(\tau_2^0))^2=\alpha \beta \gamma.
\end{equation}

By the Gross-Zagier singular moduli formula (\cite{GZSingular}, \cite{Dorman}, or  \cite[Theorem 1.1]{YYCM}), we have $\ord_2(\norm(j(\tau_1^0)-j(\tau_2^0)))=0$.  On the other hand,
 Theorem \ref{maintheorem}  implies
$$\ord_2(\alpha)=
\left\{\begin{array}{ll}
8h_1h_2&\text{if }~d_2\equiv 1\pmod 8,\\
2h_1h_2&\text{if } ~d_2\equiv 0\pmod 4.
\end{array}\right.$$
Now Table \ref{table6} follows from  Table \ref{table5} and (\ref{eq:NormRelation}).
\end{proof}

{\bf Proof of Theorem \ref{theo:Difference}:} Now we are ready to prove Theorem \ref{theo:Difference}. If one of the $d_i \equiv 5 \pmod 8$, for example, we assume $d_2 \equiv 5 \pmod 8$. Then Lemma \ref{lem:norm} and Proposition \ref{prop:norm} imply
$$\mathcal N(d_1,d_2)=\left\{\alpha,~\frac{\alpha}{\norm(\lambda(\tau_1^0))^{a(d_2)}},~\frac{\alpha}{\norm(1-\lambda(\tau_1^0))^{a(d_2)}}\right\}.$$
Proposition \ref{prop:norm}  shows  that $|\mathcal N(d_1,d_2)|=1$ or $2$ depending on whether $d_1 \equiv 5 \pmod 8$ or not.

Now assume that $d_i \not\equiv 5 \pmod 8$.  Then we may assume $d_1 \equiv 1 \pmod  8$  as $(d_1, d_2)=1$. Lemmas  \ref{lem:norm}, \ref{lem:numerator}, and Proposition \ref{prop:norm}  imply
\begin{align*}
\mathcal N(d_1,d_2)&=\left\{\alpha,~\frac{2^{-8h_1h_2}\alpha}{\norm(\lambda(\tau_2^0))^{2h_1}},~\frac{2^{-8h_1h_2}\alpha}{ \norm(1-\lambda(\tau_2^0))^{2h_1}}, ~
\tilde\alpha,~\frac{2^{-8h_1h_2}\tilde\alpha}{\norm(\lambda_i(\tau_2^0))^{2h_1}},~\frac{2^{-8h_1h_2}\tilde\alpha}{\norm(1-\lambda_i(\tau_2^0))^{2h_1)}}
\right\},
\end{align*}
where $\tilde\alpha = \norm(\lambda(\tau_1^0) -\lambda_i(\tau_2^0)) \ne \alpha $ for some $i$ as in Lemma \ref{lem:numerator}.

When $d_2 \equiv 1 \pmod 8$, then Proposition \ref{prop:norm} and Lemma \ref{lem:numerator} show that
$$
\mathcal N(d_1,d_2)=\left\{\alpha,\, 2^{-16 h_1 h_2} \alpha,\, \tilde\alpha, \,  2^{-8 h_1h_2} \tilde\alpha \right\}
$$
has cardinality $4$ as their $2$-adic valuations are distinct:  $8h_1h_2, -8h_1 h_2,  -4 h_1h_2, -12h_1h_2$.

Similarly one has for $d_2 \equiv 4, 0 \pmod 8$,
$$
\mathcal N(d_1,d_2)=\left\{\alpha,\,2^{-8h_1h_2}\alpha, \,  2^{-12 h_1 h_2} \alpha,\,  \tilde\alpha, \, \, 2^{-4 h_1h_2} \tilde\alpha \right\}
$$
which has cardinality $5$.  This proves the theorem.

\section{Appendix A:  Explicit formulas for Whittaker functions} \label{sect5}

Let $F$ be a finite field extension of $\Q_p$ with ring of integers $\OO_F$ and a  uniformizer $\pi=\pi_F$. We denote by $|\cdot|$ the valuation of $F$. Let $E$ be a quadratic \'etale extension of $F$, including $F \times F$, with ring of integers $\OO_E$ and a uniformizer $\pi_E$, and let $\chi=\chi_{E/F}$ be the associated quadratic character of $F^\times$. Let $\psi$ be an unramified additive character of $F$, and let $\psi_E =\psi \circ \tr_{E/F}$.  Fix $0\ne \beta \in \OO_F$, and let $W=W_\beta=E$ with $F$-quadratic form $Q_\beta(x) =\beta x \bar x$. Let $\omega=\omega_{\beta, \psi}$ be the  Weil representation of $G=\SL_2(F)$ on $S(W)=S(E)$ associated to the reductive dual pair $(O(W), \SL_2)$, and let
\begin{equation}
\lambda: S(E) \rightarrow I(0, \chi) , \quad \lambda(\phi) (g) =\omega(g)\phi(0).
\end{equation}
Let $\Phi_\phi \in  I(s, \chi)$ be the associated standard section, given by
$$
\Phi_\phi(g, s) = \lambda(\phi)(g) |a(g)|^s,
$$
where $|a(g)|=|a| $ if $g =n(b) m(a) k$ with $b \in F$, $a \in F^\times$, and $k \in \SL_2(\OO_F)$.  Recall that
$$
n(b) =\kzxz {1} {b} {0} {1}, \quad m(a)=\kzxz {a} {0} {0} {a^{-1}} \quad \hbox{ and }\ w=\kzxz {0} {-1} {1} {0}.
$$
The Whittaker function  of $\Phi_\phi$ (and  thus $\phi)$ is  defined to be
\begin{equation}
W_t(g, s, \phi) =\int_{F} \Phi_\phi(w n(b)g, s) \psi(-tb) db,
\end{equation}
where $db$ is the self-dual Haar measure on $F$ with respect to $\psi$, i.e., $\vol(\OO_F)=1$. The purpose of this appendix is to compute the Whittaker function
$$
W_t(s, \phi) =W_t(1, s, \phi)
$$
explicitly in the general case, which is of  independent interest in addition to the applications in Section \ref{sect:BigCM}. Various special cases have been appeared in  \cite{Yang05}, \cite[Section 4.6]{HYbook}, and \cite{KY10}.  Notice that
$$
S(E) =\bigcup_{L \subset \OO_E,\ \hbox{\tiny ideal }} S_L,
$$
where $S_L$ is  defined  in Section \ref{sect:Borcherds}.
We can and will assume $\phi=\phi_\mu =\cha(\mu+L)$ in the following. Via shifting, we can  and will further assume
\begin{equation}
L=\OO_E,  \quad Q_\beta(x) =\beta x \bar x.
\end{equation}
So
$$
L'=\{ x\in E:\,  \psi((x, y))=1,  \hbox{ for all } y \in  \OO_E\} = (\beta\partial_{E/F})^{-1}\OO_E.
$$
We denote $o(x) =\hbox{ord}_\pi x$ for $x \in F$, and
\begin{equation} \label{eq5.4}
 o_E(\mu) = \begin{cases}
 \hbox{ord}_{\pi_E}(\mu) &\ff E/F \hbox{ non-split},
 \\
  \min(o(\mu_1), o(\mu_2)) &\ff E=F \times F,\,  \mu=(\mu_1, \mu_2).
  \end{cases}
\end{equation}
We also write $f= \ord_{\pi_E} (\partial_{E/F})$ and  $e=e(E/F)=1,\hbox{ or } 2$ for the ramification index. Define  the generalized Gauss integral  for $\mu \in E/\OO_E$ and $b \in F^\times$ to be
\begin{equation}
I_\mu(b)=\int_{\mu +\OO_E} \psi(b x \bar x) d x,
\end{equation}
where $dx$ is the standard Haar measure on $E$ with $\vol(\OO_E, dx) =1$. Let
$d_\beta x$ be the self-dual Haar measure on $W_\beta=E$ with respect to $\psi ((x, y))=\psi(\tr_{E/F} \beta x\bar y)$. Then $d_\beta x = |\beta| |d_{E/F}|^{\frac{1}2} dx$, i.e.
$\vol(\OO_E, d_\beta x) =|\beta| |d_{E/F}|^{\frac{1}2}$.

\begin{proposition} \label{prop:WhittakerFunction} Let $\psi$ be an unramified additive character of $F$, and $W_\beta =E$ with $F$-quadratic form $Q_\beta(z) =\beta z \bar z$ with $0 \ne \beta \in \OO_F$. Let $L=\cha(\OO_E)$ and $\mu \in  L'$.  Then
\begin{equation} \label{eqWSummation}
\frac{W_t(s, \phi_\mu)}{\gamma(\beta) \vol(\OO_E, d_\beta x)}
 = \sum_{n=0}^\infty J_\mu(n) q^{n(1-s)},
\end{equation}
where $\gamma(\beta)=\gamma(W_\beta)$ is the local Weil index of $W_\beta$ and
\begin{equation} \label{eqJn}
J_\mu(n) =\begin{cases}
 \int_{\OO_F} I_\mu(b\beta) \psi(-tb) db    &\ff n =0,
 \\
  \int_{\OO_F^\times} I_\mu(b\beta \pi^{-n}) \psi(-tb\pi^{-n} ) db  &\ff n >0.
\end{cases}
\end{equation}
Moreover,
$$
J_\mu(0) = \cha(\OO_F) (Q_\beta(\mu) -t).
$$
\end{proposition}
\begin{proof} Since
$$
|a(wn(b))|=\begin{cases}
    1 &\ff b \in \OO_F,
    \\
     |b|^{-1} &\ff b \notin \OO_F,
     \end{cases}
$$
 one has by definition,
\begin{align*}
\frac{W_t(s, \phi_\mu)}{\gamma(\beta) \vol(\OO_E,d_\beta x) }
 &= \frac{1}{\vol(\OO_E,d_\beta x) }  \int_{F} \int_{E} \psi(bQ_\beta(x)) \phi_\mu(x) d_\beta x  |a(wn(b))|^s \psi(-t b) db
 \\
 &=\int_{F} I_\mu(b \beta) \psi(-tb) |a(wn(b))|^s db
 \\
 &=\sum_{n=0}^\infty J_\mu(n) q^{n(1-s)}.
\end{align*}
Next,
\begin{align*}
J_\mu(0)  &=\int_{\mu +\OO_E} dx  \int_{\OO_F}  \psi(b (Q_\beta(x) -t) ) db
\\
 &= \int_{\mu+\OO_E} \cha(\OO_F) (Q_\beta(x) -t) dx
 \\
 &=\cha(\OO_F)( Q_\beta(\mu) -t)
\end{align*}
as claimed.
\end{proof}

\begin{lemma}  \label{lem:GaussIntegral} Let $\mu \in E/\OO_E$ and $b \in F^\times$.
\begin{enumerate}

\item When $\mu =0$, one has
$$
I_0(b) =\begin{cases}
    1 &\ff  b \in \OO_F,
    \\
    \chi(b) |b|^{-1}  &\ff    b \notin \OO_F,\ E/F \hbox{ unramified},
    \\
    0 &\ff   -f < o(b) <0,\ E/F \hbox{ ramified},
  \\
  \chi(b) |b|^{-1} q^{-\frac{f}2} \gamma(W_1) &\ff o(b) \le -f,\ E/F \hbox{ ramified}.
  \end{cases}
$$

\item When $\mu \notin \OO_E$ and  $E/F$ is unramified (inert or split),  one has
$$
I_\mu(b)  =
    \psi(b \mu \bar\mu) \cha(\OO_E)(b\mu) .
$$

\item When  $\mu \notin \OO_E$ and  $E/F$ is ramified, one has
 $$
 I_\mu(b)
 =\begin{cases}
 \psi(b \mu \bar\mu) \cha(\partial_{E/F}^{-1})(b \mu) &\ff b \in \OO_F,
 \\
  |b|^{-1} \sum_{a \in \OO_E/\pi_E^{m}} \psi(b(\mu+a) \overline{(\mu+a)}) &\ff  0< m \le o_E(\mu) + f,
  \\
  0 &\ff m > o_E(\mu) + f.
  \end{cases}
$$
Here $m=-o(b)>0$.
\end{enumerate}
\end{lemma}
\begin{proof} We first consider the case $\mu=0$. The case $E=F\times F$ is trivial. Assume that $E$ is a quadratic field extension of $F$. Then \cite[Lemma 4.6.1]{HYbook} implies
$$
I_0(b)= \int_{\OO_E} \psi(b z \bar z) dz  =C \int_{\OO_F} \psi(bx) (1+\chi(x)) dx,
$$
with  $C=1$ or $\frac{1}{2}(1+q^{-1})$ depending on whether $E/F$ is ramified or not. When  $E/F$ is ramified, the calculation in \cite[Page 60]{HYbook} gives (1) for $\mu=0$.  When $E/F$ is  inert, $\chi$ is unramified and $\chi(\pi)=-1$. For $b \notin \OO_F$, write $m=-o(b) >0$, then one has
\begin{align*}
\int_{\OO_F} \psi(bx) \chi(x) dx
 &=\sum_{n=0}^\infty q^{-n} \chi(\pi)^{n} \int_{\OO_F^\times} \psi(\pi^{n} b x)  dx
 \\
  &=\sum_{n=0}^\infty q^{-n} \chi(\pi)^{n} (\cha(\OO_F)(\pi^n b) -q^{-1} \cha(\OO_F)(\pi^{n+1} b))
  \\
  &= \left[(1-q^{-1}) \sum_{n=m}^{\infty} (-q)^n \right]- q^{-1} (-q)^{-m+1}
  \\&=
  \frac{2}{1+q^{-1}} \chi(b) |b|^{-1}.
\end{align*}
So one has for $b \notin \OO_F$,
$$
I_0(b)= \chi(b) |b|^{-1}
$$
as claimed. For $b \in  \OO_F$, it is trivial.  This takes care  the case $\mu \in  \OO_E$.

Next, we assume that $\mu \notin \OO_E$ and $E/F$ is unramified.  When $b \mu \in \OO_E$, one has $b \in \OO_F$ and
$$
I_\mu(b) = \psi(b\mu \bar\mu) \int_{\OO_E} \psi(\tr_{E/F}( b\bar\mu z)) \psi(b z \bar z) dz = \psi(b\mu \bar\mu)\int_{\OO_E} \psi(b z \bar z) dz =\psi(b \mu\bar\mu) .
$$
When $b \mu \notin \OO_E$, and $b \in \OO_F$, the same computation gives
$$
I_\mu(b) = \psi(b\mu \bar\mu) \int_{\OO_E} \psi(\tr_{E/F}( b\bar\mu z)) ) dz=0
$$
Finally, when $b \mu \notin \OO_E$, and $b \notin \OO_F$, A substitution $z \mapsto z+ b^{-1} u$, $u \in \OO_E^\times$, gives
$$
I_\mu(b) = \psi(\tr_{E/F}(u\bar\mu)) I_\mu(b).
$$
Since there is $u \in \OO_F^\times$ such that $\psi(\tr_{E/F}(u\mu)) \ne 1$, one sees that $I_\mu(b)=0$. This proves (2).

Finally, assume that $\mu \notin \OO_E$ and $E/F$ is ramified. The case $b \in \OO_F$ is the same as in the proof of (2). Now assume $b \notin \OO_F$, so $m=-o(b) >0$.
\begin{align*}
I_\mu(b)
 &= \sum_{a \in \OO_E/\pi_E^{m}}
   \int_{\mu +a + \pi_E^{m} \OO_E} \psi(b x \bar x) d x
 \\
  &=  \sum_{a \in \OO_E/\pi_E^{m}} q_E^{-m}  \psi(b(\mu+a) \overline{(\mu+a)})
     \int_{\OO_E} \psi_E( b \pi_E^{m} \overline{(\mu+a)} x) d x
 \\
  &= q_E^{-m}\sum_{a \in \OO_E/\pi_E^{m}}  \psi(b(\mu+a) \overline{(\mu+a)})\cha(\partial_{E/F}^{-1}) (b \pi_E^{m} \overline{(\mu+a)})
  \\
  &= q^{-m}\sum_{a \in \OO_E/\pi_E^{m}}  \psi(b(\mu+a) \overline{(\mu+a)})\cha(\partial_{E/F}^{-1}) (b \pi_E^{m} \bar\mu).
\end{align*}
Now (3) is clear.

\end{proof}

We denote
\begin{equation}
\alpha(\mu, t) =Q_\beta(\mu) -t= \beta \mu \bar\mu -t, \quad a(\mu, t) =o(Q_\beta(\mu) -t)= o(\beta \mu \bar\mu -t).
\end{equation}

\begin{proposition} \label{prop7.3} Assume that $E/F$ is unramified and $0 \ne \mu \in E/\OO_E$. Then $W_t(s, \phi_\mu)=0$ unless $\alpha(\mu, t)\in  \OO_F$, i.e., $a(\mu, t) \ge 0$. Moreover,
$$
\frac{W_t(s, \phi_\mu)}{\gamma(\beta) \vol(\OO_E, d_\beta x)}
=\begin{cases}
  (1-q^{-s}) \sum_{0 \le n \le a(\mu, t)} q^{n(1-s)}  &\ff  0 \le a(\mu, t) <  o_E(\beta \mu),
  \\
  (1-q^{-s}) \sum_{0 \le n < o_E(\beta \mu)} q^{n(1-s)} + q^{o_E(\beta \mu) (1-s) } &\ff a(\mu, t) \ge  o_E(\beta \mu).
\end{cases}
$$
In  particular,  $W_t(0, \phi_\mu)=0$ if and only if $ a(\mu, t) <  o_E(\beta \mu)$. In such a case,
$$
\frac{W_t'(0, \phi_\mu)}{\gamma(\beta) \vol(\OO_E, d_\beta x)}
=\begin{cases}
 (\sum_{ 0 \le n \le a(\mu, t)} q^n)  \log q  &\ff 0 \le a(\mu, t) <  o_E(\beta \mu),
  \\
    0 &\ff a(\mu, t) < 0.
    \end{cases}
$$
\end{proposition}
\begin{proof} Since $E/F$ is unramified,  one has $\partial_{E/F} =\OO_E$, $f=0$ and $\pi_E=\pi$.  When $ 0 < n  \le o(\beta)$,  one has for $b \in \OO_F^\times$ by Lemma \ref{lem:GaussIntegral}
\begin{align*}
I_\mu(\beta b \pi^{-n})&=\psi(\pi^{-n} b Q_\beta(\mu)) \cha(\OO_E) (b  \mu  \beta \pi^{-n}).
\end{align*}
So $J_\mu(n) =0$ for $ n > o_E(\beta \mu)$.
When $0 < n \le o_E(\beta \mu)$, one has
\begin{align*}
J_\mu(n) &=\int_{\OO_F^\times} \psi(\pi^{-n} b (Q_\beta(\mu) -t)) db
\\
 &= \cha(\pi^n\OO_F) (Q_\beta(\mu) -t)) - q^{-1} \cha(\pi^{n-1}\OO_F) (Q_\beta(\mu) -t)).
\end{align*}
  So
\begin{align*}
&\frac{W_t(s, \phi_\mu)}{\gamma(\beta) \vol(\OO_E, d_\beta x)}
 =\cha(\OO_F) (Q_\beta(\mu) -t))    +
   \\
    &\qquad \sum_{0 < n \le o_E(\beta \mu) } \left[\cha(\pi^n\OO_F) (Q_\beta(\mu) -t)) - q^{-1} \cha(\pi^{n-1}\OO_F) (Q_\beta(\mu) -t))\right] q^{n(1-s)}.
\end{align*}
Now  a simple  calculation proves  the proposition.
\end{proof}

\begin{proposition} \label{prop7.4}
Assume that $E/F$ is ramified with $f =o_E(\partial_{E/F})=o(d_{E/F})$ and $0\ne \mu \in  E/\OO_E$.  Then $W_t(s, \phi_\mu) =0$ unless $a(\mu, t)  \ge 0$.  Write $a= o(\beta) + [\frac{1}2(o_E(\mu) +f)]  $, where $[x]$ is the largest integer less than or equal to $x$.
\begin{enumerate}
\item  \quad When $o_E(\mu) +f \le 0$, i.e.,  $\mu\notin  \pi_E \partial_{E/F}^{-1}$,   one has
$$
\frac{W_t(s, \phi_\mu)}{\gamma(\beta) \vol(\OO_E, d_\beta x)}
=\begin{cases}
 (1-q^{-s}) \sum_{0 \le n \le a(\mu, t) } q^{n(1-s)}  &\ff a(\mu, t) < a,
 \\
  (1-q^{-s}) \sum_{0 \le n < a } q^{n(1-s)}  + q^{a(1-s)}  &\ff a(\mu, t)  \ge  a.
 \end{cases}
 $$

 \item \quad When  $ o_E(\mu) + f >0$, i.e.,  $\mu \in  \pi_E\partial_{E/F}^{-1}-\OO_E$,   one has
 $$
 \frac{W_t(s, \phi_\mu)}{\gamma(\beta) \vol(\OO_E, d_\beta x)} = (1-q^{-s}) \sum_{0 \le n \le a(\mu, t) } q^{n(1-s)}
 $$
 if $a(\mu, t) < o(\beta)$ and
\begin{align*}
\frac{W_t(s, \phi_\mu)}{\gamma(\beta) \vol(\OO_E, d_\beta x)}
 &=(1-q^{-s})  \sum_{0 \le n < o(\beta) } q^{n(1-s)}  + q^{o(\beta)(1-s)}
 \\
  &\qquad + \sum_{o(\beta) < n \le o(\beta) + o_E(\mu) +f} q^{o(\beta)} q^{-ns}
    \sum_{a \in  \OO_E/\pi_E^{n-o(\beta)}}
    \\
    &\qquad
     \left[\cha(\pi^n\OO_F) (\alpha(\mu +a, t) ) - q^{-1} \cha(\pi^{n-1}\OO_F) (\alpha(\mu +a, t))\right]
\end{align*}
if $a(\mu, t) \ge o(\beta)$.
\item  One has $W_t(0,\phi_\mu) =0$ if and only if $a(\mu, t) < \min(o(\beta),  a)$. In such a case,
$$
\frac{W_t'(0, \phi_\mu)}{\gamma(\beta) \vol(\OO_E, d_\beta x)}
=\begin{cases}
 (\sum_{ 0 \le n \le a(\mu, t)} q^n)  \log q  &\ff 0 \le a(\mu, t) <  \min(o(\beta),  a),
  \\
    0 &\ff a(\mu, t) < 0.
    \end{cases}
$$

\end{enumerate}

\end{proposition}

\begin{remark} When $p\ne 2$, $f=1$, the second sum in Proposition \ref{prop7.4}(2) does not appear, the resulting formula is the same as the unramified case.  In general, $f \le 1+ o(2)$, the last sum is manageable when $o(2) =\ord_\pi 2$ is small. However, the ramification index $f$   can be arbitrarily large when $F/\Q_2$ has arbitrary large ramification, and the  sum can  become complicated. We refer to \cite[Appendix A]{Yang04} for an explicit way to compute $f$.
\end{remark}

\begin{proof} Assume $o_E(\mu) +f \le 0$ first. For $b \in \OO_F^\times$,  Lemma \ref{lem:GaussIntegral} implies
$$
I_\mu(\beta b\pi^{-n})=
\begin{cases}
  \psi(Q_\beta(\mu) b \pi^{-n} ) &\ff n \le a,
  \\
  0 &\ff n > a.
  \end{cases}
$$
When  $n \le a$, one has then
\begin{align*}
J_\mu(\beta b \pi^{-n})
& =\int_{\OO_F^\times} \psi( \pi^{-n}(Q_\beta(\mu) -t) b) db
\\
&=\begin{cases}
 1-q^{-1} &\ff n \le \min (a(\mu, t) , a),
 \\
 -q^{-1}  &\ff n =a(\mu, t) +1 \le a,
 \\
  0  &\ff n > a(\mu, t) +1.
 \end{cases}
\end{align*}
Applying Proposition \ref{prop:WhittakerFunction}, we obtain  (1).

Now assume that $o_E(\mu) + f >0$. When $ m=n-o(\beta) \le 0$, we have with the same calculation
$$
J_\mu(n) =  \begin{cases}
  1- q^{-1}  &\ff n \le \min(a(\mu, t), o(\beta)),
  \\
   -q^{-1} &\ff n=a(\mu, t) +1 \le o(\beta),
   \\
   0 &\ff n > a(\mu, t) +1.
   \end{cases}
$$
When $m >  0$, i.e, $n \ge o(\beta)$, Lemma \ref{lem:GaussIntegral} implies
$$
I_\mu(\beta b \pi^{-n})
 =q^{-m} \sum_{a \in \OO_E/\pi_E^{n-o(\beta)}} \psi(\pi^{-n} b Q_\beta(\mu+a)),
$$
and so
$$
J_\mu(n) = q^{-m}
\sum_{a \in \OO_E/\pi_E^{n-o(\beta)}}\left[ \cha(\pi^n\OO_F) (Q_\beta(\mu+a) -t)) - q^{-1} \cha(\pi^{n-1}\OO_F) (Q_\beta(\mu+a) -t))   \right].
$$
Now (2) is clear.   Claim (3) follows from (1) and (2).

\end{proof}

\begin{corollary}  \label{cor:ConstantTerm} Assume $\mu \in E-\OO_E$. Then  $W_0(s, \phi_\mu) =0$ unless $a(\mu, 0) =o(\beta \mu \bar\mu) \ge 0$. When  $a(\mu, 0) \ge 0$,  the following is true.
\begin{enumerate}
\item When  $E/F$ is non-split, one has
$$
\frac{W_0(s, \phi_\mu)}{\gamma(\beta) \vol(\OO_E, d_\beta x)} = (1-q^{-s}) \sum_{0 \le n \le a(\mu, 0) } q^{n(1-s)}.
$$
In  particular,  $W_0(0, \phi_\mu) =0$, and
$$
\frac{W_0'(0, \phi_\mu)}{\gamma(\beta) \vol(\OO_E, d_\beta x)}= \sum_{0 \le n \le a(\mu, 0) } q^{n} \log q.
$$

\item  When  $E=F \oplus F$ is split over $F$ and $\mu=(\mu_1, \mu_2) \in (F^\times)^2$, one has
$$
\frac{W_0(s, \phi_\mu)}{\gamma(\beta) \vol(\OO_E, d_\beta x)}
= \begin{cases}
   (1-q^{-s}) \sum_{0 \le n \le a(\mu, 0) } q^{n(1-s)} &\ff \mu_1, \mu_2 \notin \OO_F,
   \\
    (1-q^{-s}) \sum_{0 \le n < o_E(\beta \mu) } q^{n(1-s)} + q^{o_E(\beta\mu)(1-s)} &\ff \text{only one of $\mu_1$ and $\mu_2$ is in $\OO_F$}.
\end{cases}
$$
In particular
$$
\frac{W_0'(0, \phi_\mu)}{\gamma(\beta) \vol(\OO_E, d_\beta x)}
 =\begin{cases}
  \sum_{0 \le n \le a(\mu, 0) } q^{n} \log q  &\ff \mu_1, \mu_2 \notin \OO_F,
  \\
   \sum_{0 \le n < o_E(\beta \mu) } q^{n} \log q  - q^{o_E(\beta \mu)} \log q &\ff \text{only one of $\mu_1$ and $\mu_2$ is in $\OO_F$}.
\end{cases}
$$

\end{enumerate}

\end{corollary}
\begin{proof} Since $\mu \notin \OO_E$,   $a(\mu, 0) = o(\beta \mu \bar\mu) < o_E(\beta \mu)$ when $E/F$ is inert, and $a(\mu, 0) < \min(o(\beta), a)$ when $E/F$ is ramified.  This implies (1) by Propositions \ref{prop7.3} and \ref{prop7.4}.  The first part of (2) is the same.  On the other hand,  if one of $\mu_i$, say,  $\mu_1 \in \OO_F$, one has
$$
a(\mu, 0) = o(\beta \mu_1 \mu_2) \ge o(\beta \mu_2) = o_E(\beta \mu).
$$
Proposition \ref{prop7.3} implies this case too.

\end{proof}

The case $\mu=0$  is already computed in \cite[Section 4.6]{HYbook} using a slightly different method. There are some  minor errors there. We record the correct results here for the convenience of the reader.

\begin{proposition} \label{prop7.6}
(\cite[Proposition 4.6.2]{HYbook}) Assume that $E/F$ is unramified and let  $ t \in F^\times$. Then $W_t(s, \phi_0) =0$ unless $t \in  \OO_F$. Assume $t \in \OO_F$. Recall that $\chi=\chi_{E/F}$ is the quadratic character of $F^\times$ associated to $E/F$.

\begin{enumerate}
\item \quad When  $o(t) < o(\beta)$, one has
$$
\frac{W_t(s, \phi_0)}{\gamma(\beta) \vol(\OO_E, d_\beta x)}
 = (1-q^{-s}) \sum_{0 \le n \le o(t)} q^{n(1-s)}.
$$
In particular, $W_t(0, \phi_0) =0$.

\item \quad When $o(t) \ge o(\beta)$,  and $t \ne 0$, one has
$$
\frac{W_t(s, \phi_0)}{\gamma(\beta) \vol(\OO_E, d_\beta x)}=(1-q^{-s}) \sum_{0 \le n \le o(\beta)-1} q^{n(1-s)}
 + \frac{q^{o(\beta)(1-s)}}{L(s+1,\chi)}\sum_{0 \le n \le o(t/\beta)} (\chi(\pi) q^{-s})^n.
$$
In particular,
$$
L(1,\chi)\frac{W_t(0, \phi_0)}{\gamma(\beta) \vol(\OO_E, d_\beta x)}
=\begin{cases}
 q^{o(\beta)}(o(t/\beta) +1) &\ff E/F \hbox{ is split},
 \\
 q^{o(\beta)}\frac{1+(-1)^{o(t/\beta)}}2 &\ff  E/F \hbox{ is inert}.
 \end{cases}
$$
Moreover,  $W_t(0, \phi_0) =0$ if and only if $E/F$ is  inert and $o(t/\beta)$ is  odd.  When this is the case, one has
\[L(1,\chi)\frac{W'_t(0, \phi_0)}{\gamma(\beta) \vol(\OO_E, d_\beta x)}
=\left[q^{o(\beta)}\frac{\ord_\pi(t/\beta)+1}{2}+\frac{1-q^{\ord_\pi(\beta)}}{q(1-q^{-2})}\right]\log(q).\]

\item \quad When $t=0$, one has
$$
\frac{W_0(s, \phi_0)}{\gamma(\beta) \vol(\OO_E, d_\beta x)}=(1-q^{-s}) \sum_{0 \le n \le o(\beta)-1} q^{n(1-s)}+q^{o(\beta)(1-s)}\frac{ L(s, \chi)}{L(s+1,\chi)}.
$$
\begin{proof}
First assume $t\neq 0$. By Proposition \ref{prop:WhittakerFunction} and Lemma \ref{lem:GaussIntegral}(2), $J_0(0)=1$ and for $0<n\leq o(\beta)$,
\begin{align*}
J_0(n) &=\int_{\OO_F^\times} \psi(-\pi^{-n} bt) db
\\
 &=\begin{cases}\cha(\pi^n\OO_F) (t) - q^{-1} \cha(\pi^{n-1}\OO_F) (t)&\ff t\neq 0,\\
                             1-q^{-1}&\ff t=0.
\end{cases}
\end{align*}
and for $n>o(\beta)$,
\begin{align*}
J_0(n) &=\int_{\OO_F^\times}\chi(b\beta\pi^{-n})|b\beta\pi^{-n}|^{-1} \psi(-\pi^{-n} bt) db
\\
 &=\begin{cases}\chi(\pi)^{o(\beta)-n}q^{o(\beta)-n}(\cha(\pi^n\OO_F) (t) - q^{-1} \cha(\pi^{n-1}\OO_F) (t))&\ff t\neq 0,\\ \chi(\pi)^{o(\beta)-n}q^{o(\beta)-n}(1-q^{-1})&\ff t=0.
\end{cases}
\end{align*}
When  $t\neq 0$ and $o(t)<o(\beta)$ we have  by (\ref{eqWSummation}),
\begin{eqnarray*}
\frac{W_t(s, \phi_0)}{\gamma(\beta) \vol(\OO_E, d_\beta x)}&=&1+\sum_{0<n\leq o(\beta)}(1-q^{-1})q^{n(1-s)}-q^{-1}q^{(o(t)+1)(1-s)}\\
&=&(1-q^{-s})\sum_{0 \le n \le o(t)} q^{n(1-s)}.
\end{eqnarray*}
When f $t\neq 0$ and $o(t)\geq o(\beta)$, we have 
\begin{eqnarray*}
\frac{W_t(s, \phi_0)}{\gamma(\beta) \vol(\OO_E, d_\beta x)}&=&1+\sum_{0<n\leq o(\beta)}(1-q^{-1})q^{n(1-s)}\\&&+\sum_{o(\beta)<n\leq o(t)}\chi(\pi)^{o(\beta)-n}q^{o(\beta)-n}(1-q^{-1})q^{n(1-s)}\\&&-q^{-1}\chi(\pi)^{o(\beta)-(o(t)+1)}q^{o(\beta)-(o(t)+1)}q^{(o(t)+1)(1-s)}\\
&=&(1-q^{-s}) \sum_{0 \le n \le o(\beta)-1} q^{n(1-s)}
 + \frac{q^{o(\beta)(1-s)}}{L(s+1,\chi)}\sum_{0 \le n \le o(t/\beta)} (\chi(\pi) q^{-s})^n.
\end{eqnarray*}
When  $t=0$, we have 
\begin{eqnarray*}
\frac{W_t(s, \phi_0)}{\gamma(\beta) \vol(\OO_E, d_\beta x)}&=&1+\sum_{0<n\leq o(\beta)}(1-q^{-1})q^{n(1-s)}\\&&+\sum_{o(\beta)<n}\chi(\pi)^{o(\beta)-n}q^{o(\beta)-n}(1-q^{-1})q^{n(1-s)}\\
&=&(1-q^{-s}) \sum_{0 \le n \le o(\beta)-1} q^{n(1-s)}
 +\frac{q^{o(\beta)(1-s)} L(s, \chi)}{L(s+1,\chi)} .
\end{eqnarray*}

\end{proof}

\end{enumerate}
\end{proposition}

\begin{proposition} \label{prop7.7}
(\cite[Proposition 4.6.3]{HYbook}) Assume that $E/F$ is ramified and let  $ t \in F$. Then $W_t(s, \phi_0) =0$ unless $t \in  \OO_F$. Assume $t \in \OO_F$. Recall that $\chi=\chi_{E/F}$ is the quadratic character of $F^\times$ associated to $E/F$.

\begin{enumerate}
\item \quad When  $o(t) < o(\beta)$, one has
$$
\frac{W_t(s, \phi_0)}{\gamma(\beta) \vol(\OO_E, d_\beta x)}
 = (1-q^{-s}) \sum_{0 \le n \le o(t)} q^{n(1-s)}.
$$
In particular, $W_t(0, \phi_0) =0$.

\item \quad When $o(t) \ge o(\beta)$,  and $t \ne 0$, one has
$$
\frac{W_t(s, \phi_0)}{\gamma(\beta) \vol(\OO_E, d_\beta x)}=|\beta|^{s-1} \left(1+\chi(t/\beta) q^{-(o(t/\beta) +f)s}\right) + (1-q^{-s}) \sum_{0 \le n \le o(\beta)-1} q^{n(1-s)}.
$$
In particular,
$$
\frac{W_t(0, \phi_0)}{\gamma(\beta) \vol(\OO_E, d_\beta x)}
=|\beta|^{-1} (1+\chi(t/\beta)).
$$
Moreover,  $W_t(0, \phi_0) =0$ if and only if $\chi(t/\beta)=-1$, i.e., $W_\beta$ does not represent $t$.  When this is the case,
\[L(1,\chi)\frac{W'_t(0, \phi_0)}{\gamma(\beta) \vol(\OO_E, d_\beta x)}
=\left[\ord_\pi(t/\beta)+f+\frac{1-q^{\ord_\pi(\beta)}}{q(1-q^{-1})}\right]\log(q).\]

\item \quad When $t=0$, one has
$$
\frac{W_0(s, \phi_0)}{\gamma(\beta) \vol(\OO_E, d_\beta x)}= |\beta|^{s-1}+(1-q^{-s}) \sum_{0 \le n \le o(\beta)-1} q^{n(1-s)}.
$$
\end{enumerate}
\end{proposition}

\subsection{Variants} The first variant is  the case  when    $\psi$ has conductor $n(\psi)=n \ne 0$, where
$$
n(\psi) = \{ n \in  \Z:\,  \psi|_{\pi^n \OO_F} =1\}.
$$
Let $\psi'(x) = \psi(a x)$ for some $a \in F$ with $o(a) =n$, then $n(\psi')=0$, i.e., $\psi'$ is unramified. Let $\beta'= a^{-1}  \beta$, then  the Weil representation $\omega_{W_\beta, \psi} = \omega_{W_{\beta'}, \psi'}$, and $d_{\psi'} b = |a|^{\frac{1}2} d_\psi b$. So
\begin{align*}
W_{t}^\psi (s,  \phi)
 &= \int_F \omega_{W_\beta, \psi} (w n(b) )(0) |a(wn(b))|^s \psi(-tb) d_\psi b
 \\
  &=|a|^{-\frac{1}2} \int_F \omega_{W_\beta', \psi'} (w n(b) )(0) |a(wn(b))|^s \psi¡®(-a^{-1} tb) d_{\psi'}  b
  \\
  &=|a|^{-\frac{1}2}W_{\frac{t}a}^{\psi'}(s, \phi).
  \end{align*}
  Therefore, we have for every  $\phi \in S(E)$
  \begin{equation} \label{eq:ChangPsi}
  W_{t}^\psi (s,  \phi) = |a|^{-\frac{1}2}W_{\frac{t}a}^{\psi'}(s, \phi),
  \end{equation}
where the right hand side is with respect to $(W_{\beta'}, \psi')$.

The second variant is  the case when $L=\pi_E^c \OO_E$ with quadratic form $Q_\beta(x) =\beta x \bar x$  (assuming $\psi$ to be unramified). In such a case, let $\tilde L = \OO_E$ with quadratic form $Q_{\tilde\beta}(x) = \tilde\beta x \bar x$ and $ \tilde\beta = \pi_E^{2c} \beta$. Then
$$
(L, Q_\beta) \cong (\tilde L, Q_{\tilde\beta}),
$$
and
$$
W_{t}(s, \phi_\mu) = W_{t}(s, \phi_{\tilde\mu})
$$
with
$\phi_\mu=\cha(\mu+L)$ and $\phi_{\tilde\mu} =\cha( \tilde\mu+\OO_E)$.

\section{Appendix B:  Computation Data}

Theorem  \ref{maintheorem} can easily be used to calculate the norm $\norm(\lambda(\frac{d_1+\sqrt{d_1}}2) - \lambda(\frac{d_2+\sqrt{d_2}}2))$ up to sign $\pm 1$. Here are some examples with  $\tau_i= \frac{d_i +\sqrt{d_i}}2$.

\begin{table}[h]
\caption{$d_1\equiv 5\pmod 8,~d_2\equiv 5\pmod 8$}
\begin{tabular}{|c|c|c||c|c|c|}
\hline \rule[-3mm]{0mm}{8mm}
$d_1 $ & $d_2 $ &$\mathrm{N}(\lambda(\tau_1)-\lambda(\tau_2))$ &$d_1 $ & $d_2 $ &$\mathrm{N}(\lambda(\tau_1)-\lambda(\tau_2))$ \\
\hline $-3$ &$-11$ &$2^{14}$      &$-11$ &$-19$  &$2^{48} 13^6$\\
\hline $-3$ &$-19$ &$2^{14}3^6$   &$-11$ &$-35$  &$2^{96} 7^{12}19^6$\\
\hline $-3$ &$-35$ &$2^{28}5^6$   &$-19$ &$-35$  &$2^{96} 19^6 31^6 41^6$\\
\hline
\end{tabular}
\end{table}

\begin{table}[h] \label{table8}
\caption{$d_1\equiv 5\pmod 8,~d_2\equiv 0,4\pmod 8$}
\begin{tabular}{|c|c|c||c|c|c|}
\hline \rule[-3mm]{0mm}{8mm}
$d_1 $ & $d_2 $ &$\mathrm{N}(\lambda(\tau_1)-\lambda(\tau_2))$ &$d_1 $ & $d_2 $ &$\mathrm{N}(\lambda(\tau_1)-\lambda(\tau_2))$ \\
\hline $-3$   &$-4$   &$3$          &$-19$ &$-4$  &$2^{48} 13^6$\\
\hline $-3$   &$-8$   &$5^2$        &$-19$ &$-9$  &$2^{96} 7^{12}19^6$\\
\hline $-3$   &$-20$  &$5^2 11^2$   &$-19$ &$-20$  &$2^{96} 19^6 31^6 41^6$\\
\hline $-11$  &$-4$   &$7^2 11$     &$-19$ &$-24$  &$3^{12} 13^4 89^2 113^2$\\
\hline $-11$  &$-8$   &$7^4 13^2$   &$-35$ &$-4$  &$7^2 19^2 31^2$\\
\hline $-11$  &$-20$  &$11^2 13^4 17^4 19^2$  &$-35$ &$-8$  &$5^6 7^4 23^4 61^2$\\
\hline $-11$  &$-24$  &$13^4 17^2 19^4 41^2$  &$-35$ &$-24$  &$19^4 23^4 37^4 41^2 43^4 67^4 89^2$\\
\hline
\end{tabular}
\end{table}

\begin{table}[h]
\caption{$d_1\equiv 1\pmod 8,~d_2\equiv 5\pmod 8$}
\begin{tabular}{|c|c|c||c|c|c|}
\hline \rule[-3mm]{0mm}{8mm}
$d_1 $ & $d_2 $ &$\mathrm{N}(\lambda(\tau_1)-\lambda(\tau_2))$ &$d_1 $ & $d_2 $ &$\mathrm{N}(\lambda(\tau_1)-\lambda(\tau_2))$ \\
\hline $-7$  &$-3$   &$3^2 5^2$             &$-23$ &$-11$ &$7^{12} 11^4 17^4 19^4 43^2 61^2$\\
\hline $-7$  &$-11$  &$7^2 13^2 17^2 19^2$  &$-23$ &$-19$ &$19^2 37^2 53^2 67^2 79^2 89^2 97^2 103^2 107^2 109^2$\\
\hline $-7$  &$-19$  &$3^{14} 13^2 31^2$    &$-23$ &$-35$ &$5^{18} 7^{12} 19^8 23^4 37^4 43^4 53^4 67^4 181^2 199^2$\\
\hline $-15$  &$-11$ &$7^8 11^2 13^4 29^2 41^2$         &$-31$ &$-3$  &$3^6 11^2 17^2 23^2$\\
\hline $-15$  &$-19$ &$3^{12} 13^4 29^2 41^2 59^2 71^2$ &$-31$ &$-11$  &$11^4 13^6 17^4 29^2 43^2 73^2 79^2 83^2$\\
\hline $-23$  &$-3$  &$5^6 11^2 17^2$       &$-31$ &$-19$  &$3^{40} 13^6 29^4 37^2 127^2$\\
\hline
\end{tabular}
\end{table}

\begin{table}[h]
\caption{$d_1\equiv 1\pmod 8,~d_2\equiv 0,4\pmod 8$}
\begin{tabular}{|c|c|c||c|c|c|}
\hline \rule[-3mm]{0mm}{8mm}
$d_1 $ & $d_2 $ &$\mathrm{N}(\lambda(\tau_1)-\lambda(\tau_2))$ &$d_1 $ & $d_2 $ &$\mathrm{N}(\lambda(\tau_1)-\lambda(\tau_2))$ \\
\hline $-7$  &$-4$   &$2^1 3^2$             &$-23$ &$-8$ &$2^6 5^6 7^4 23$\\
\hline $-7$  &$-8$  &$2^2 5^2 7$  &$-23$ &$-20$ &$2^{12} 5^6 11^4 17^2 19^4 53^2$\\
\hline $-7$  &$-20$  &$2^4 5^2 13^2 17^2$    &$-23$ &$-24$ &$2^{12} 17^4 19^2 23^2 37^2 61^2 67^2$\\
\hline $-7$  &$-24$ &$2^4 3^4 13^2 19^2$         &$-31$ &$-4$  &$2^3 3^6 11^2$\\
\hline $-15$  &$-4$ &$2^2 3^2 7^2$ &$-31$ &$-8$  &$2^6 13^2 23^2 29^2 31$\\
\hline $-15$  &$-8$  &$2^4 5^2 7^2 13^2$       &$-31$ &$-20$  &$2^{12} 11^4 13^4 17^2 37^2 53^2 73^2$\\
\hline $-23$  &$-4$  &$2^3 7^2 11^2$       &$-31$ &$-24$  &$2^{12} 3^{12} 13^4 17^4 43^2 61^2$\\
\hline
\end{tabular}
\end{table}

\begin{table}[h]
\caption{$d_1\equiv 1\pmod 8,~d_2\equiv 1\pmod 8$}
\begin{tabular}{|c|c|c||c|c|c|}
\hline \rule[-3mm]{0mm}{8mm}
$d_1 $ & $d_2 $ &$\mathrm{N}(\lambda(\tau_1)-\lambda(\tau_2))$ &$d_1 $ & $d_2 $ &$\mathrm{N}(\lambda(\tau_1)-\lambda(\tau_2))$ \\
\hline $-7$  &$-15$  &$2^{16} 3^4 5^2$      &$-15$ &$-23$ &$2^{48} 5^6 7^8 11^2$\\
\hline $-7$  &$-23$  &$2^{24} 5^6 7^2$      &$-15$ &$-31$ &$2^{48} 3^{12} 11^2 13^4 29^2$\\
\hline $-7$  &$-31$  &$2^{24} 3^{12} 13^2$  &$-23$ &$-31$ &$2^{72} 11^{10} 17^6 37^2 43^2$\\
\hline
\end{tabular}
\end{table}
\newpage


Table \ref{table8} and   Theorem  \ref{maintheorem}(1) suggest  that $\norm(\lambda(\frac{d_1+\sqrt{d_1}}2) - \lambda(\frac{d_2+\sqrt{d_2}}2))$ is a sixth power (up to sign), which  can be rephrased as  the following conjecture.

\begin{conjecture} Assume  $d_1 \equiv d_2  \equiv 5 \pmod 8$. For each  $t =\frac{m+\sqrt D}2 \in \OO_F$, let $\mathfrak p_t$ be the unique prime ideal of $F$ above $2$ with $2 \in \mathfrak p_t$.
$$
\sum_{ \substack{ t =\frac{m+\sqrt D}2  \in  \OO_F, \, \ord_{2}(m^2 -D) \hbox{\tiny odd}
 \\ |m| < \sqrt D, \,  m \equiv -1 \pmod 4  }}
  \rho_{E/F}(t \mathfrak p_t^{-1})
$$
is a multiple of $3$.
\end{conjecture}

\bibliographystyle{alpha}
\bibliography{reference}

\def\cprime{$'$}
\begin{thebibliography}{KRY99}

\bibitem[BKY12]{BKY12}
Jan~Hendrik Bruinier, Stephen~S. Kudla, and Tonghai Yang.
\newblock Special values of {G}reen functions at big {CM} points.
\newblock {\em Int. Math. Res. Not. IMRN}, (9):1917--1967, 2012.

\bibitem[Bor98]{Borcherds98}
R.~E. Borcherds.
\newblock Automorphic forms with singularities on {G}rassmannians.
\newblock {\em Invent. Math.}, 132(3):491--562, 1998.

\bibitem[Bru14]{Bruinier14}
Jan~Hendrik Bruinier.
\newblock On the converse theorem for {B}orcherds products.
\newblock {\em J. Algebra}, 397:315--342, 2014.

\bibitem[Dor88]{Dorman}
David~R. Dorman.
\newblock Special values of the elliptic modular function and factorization
  formulae.
\newblock {\em J. Reine Angew. Math.}, 383:207--220, 1988.

\bibitem[GZ85]{GZSingular}
Benedict~H. Gross and Don~B. Zagier.
\newblock On singular moduli.
\newblock {\em J. Reine Angew. Math.}, 355:191--220, 1985.

\bibitem[HY12]{HYbook}
Benjamin Howard and Tonghai Yang.
\newblock {\em Intersections of {H}irzebruch-{Z}agier divisors and {CM}
  cycles}, volume 2041 of {\em Lecture Notes in Mathematics}.
\newblock Springer, Heidelberg, 2012.

\bibitem[KRY99]{KRY1}
Stephen~S. Kudla, Michael Rapoport, and Tonghai Yang.
\newblock On the derivative of an {E}isenstein series of weight one.
\newblock {\em Internat. Math. Res. Notices}, (7):347--385, 1999.

\bibitem[Kud94]{KuSplit}
Stephen~S. Kudla.
\newblock Splitting metaplectic covers of dual reductive pairs.
\newblock {\em Israel J. Math.}, 87(1-3):361--401, 1994.

\bibitem[Kud97]{KuAnnals}
Stephen~S. Kudla.
\newblock Central derivatives of {E}isenstein series and height pairings.
\newblock {\em Ann. of Math. (2)}, 146(3):545--646, 1997.

\bibitem[KY10]{KY10}
Stephen~S. Kudla and Tonghai Yang.
\newblock Eisenstein series for {SL}(2).
\newblock {\em Sci. China Math.}, 53(9):2275--2316, 2010.

\bibitem[Mil]{Milne}
J.~S. Milne.
\newblock Class field theory.
\newblock {\em www.jmilne.org/math/CourseNotes/CFT.pdf}.

\bibitem[Yan04]{Yang04}
Tonghai Yang.
\newblock Local densities of 2-adic quadratic forms.
\newblock {\em J. Number Theory}, 108(2):287--345, 2004.

\bibitem[Yan05]{Yang05}
Tonghai Yang.
\newblock C{M} number fields and modular forms.
\newblock {\em Pure Appl. Math. Q.}, 1(2, part 1):305--340, 2005.

\bibitem[Yan16]{Yang16}
Tonghai Yang.
\newblock Rational structure of ${X}(n)$ over ${\Q}$ and explicit {G}alois
  action on {CM} points.
\newblock {\em Chinese Annals of Mathematics, Series B}, 37(6):821--832, Nov
  2016.

\bibitem[Yu]{YuPeng}
Peng Yu.
\newblock {CM} values of {S}iegel theta functions and {R}osenhain invariants on
  {S}iegel 3-fold ${X}_2(2)$ and applications to genus two curve cryptography.
\newblock {\em Preprint}.

\bibitem[YY]{YYCM}
Tonghai Yang and Hongbo Yin.
\newblock Difference of modular functions and their {CM} value factorization.
\newblock {\em Trans. AMS}, To appear.

\end{thebibliography}
\end{document}